\documentstyle[12pt]{article}
\newcommand{\sect}[1]{\section{#1}\setcounter{equation}{0}}

\font\mbn=msbm10 scaled \magstep1
\font\mbs=msbm7 scaled \magstep1
\font\mbss=msbm5 scaled \magstep1
\newfam\mbff
\textfont\mbff=\mbn
\scriptfont\mbff=\mbs
\scriptscriptfont\mbff=\mbss\def\mbf{\fam\mbff}
\def\Re{{\mbf R}}

\def\Z{{\mbf Z}}
\def\Co{{\mbf C}}

\def\Di{{\mbf D}}
\newtheorem{Th}{Theorem}[section]
\newtheorem{Lm}[Th]{Lemma}
\newtheorem{C}[Th]{Corollary}
\newtheorem{D}[Th]{Definition}
\newtheorem{Proposition}[Th]{Proposition}
\newtheorem{R}[Th]{Remark}
\newtheorem{E}[Th]{Example}
\author{Alexander Brudnyi\thanks{Research supported in part by NSERC.
}\\
Department of Mathematics and Statistics\\
University of Calgary, Calgary\\
Canada}
\title{Projections in the Space $H^{\infty}$ and the Corona Theorem for
Coverings of Bordered Riemann Surfaces}
\date{} 
\begin{document} 
\maketitle
\begin{abstract}
{Let $M$ be a non-compact connected Riemann surface of finite type, and 
$R\subset\subset M$ be a relatively compact domain such that 
$H_{1}(M,\Z)=H_{1}(R,\Z)$. Let $\tilde R\longrightarrow R$ be a covering. We 
study the algebra $H^{\infty}(U)$ of bounded 
holomorphic functions defined in some domains $U\subset\tilde R$. Our main 
result is a Forelli type theorem on projections in $H^{\infty}(\Di)$.  
}
\end{abstract}
%\tableofcontents
\sect{\hspace*{-1em}. Introduction.}
{\bf 1.1.} Let $X$ be a connected complex manifold and $H^{\infty}(X)$ be
the algebra of bounded holomorphic functions on $X$
with pointwise multiplication and with norm
$$
||f||=\sup_{x\in X}|f(x)|\ .
$$
Let $r:\tilde X\longrightarrow X$ be the universal covering of $X$.
The fundamental group $\pi_{1}(X)$ acts discretely on $\tilde X$ by
biholomorphic maps. By $r^{*}(H^{\infty}(X))\subset H^{\infty}(\tilde X)$ we 
denote the Banach subspace of functions invariant with respect to the action 
of $\pi_{1}(X)$. In this paper we describe a class of  manifolds
$X$ for which there is a linear continuous projector 
$P: H^{\infty}(\tilde X)\longrightarrow r^{*}(H^{\infty}(X))$ satisfying 
\begin{equation}\label{mainprop}
P(fg)=P(f)g\ \ \ {\rm for\ any}\ f\in H^{\infty}(\tilde X),\ 
g\in r^{*}(H^{\infty}(X)).
\end{equation}
For instance, according to Forelli [F], such $P$ exists in the case when $X$ 
is the interior of a compact bordered Riemann surface. 
(The universal covering of such $X$ is the open unit disk $\Di\subset\Co$.)
One of the possible applications of Forelli's theorem is to the solution of
the corona problem for $H^{\infty}(R)$ (for further results and references 
related to the corona problem we refer to Garnett [Ga1], Jones and 
Marshall [JM] and Slodkowski [S]). Extensions of Forelli's theorem to some
Riemann surfaces of Widom type were obtained by Carleson [Ca2] and
Jones and Marshall [JM]. In this paper we consider another more general
construction of $P$ satisfying (\ref{mainprop}). Let us formulate our result. 

Let $N\subset\subset M$ be a relatively compact domain (i.e. an open 
connected subset) in a connected Stein manifold $M$ such that
\begin{equation}\label{condit1}
\pi_{1}(N)\cong\pi_{1}(M).
\end{equation}
By ${\cal F}_{c}(N)$ we denote the class of unbranched coverings 
of $N$. Recall that any covering from ${\cal F}_{c}(N)$ 
corresponds to a subgroup of $\pi_{1}(N)$. Assume that 
the complex connected manifold $U$ admits a holomorphic embedding 
$i:U\hookrightarrow R$ into some $R\in {\cal F}_{c}(N)$. 
Let $i_{*}:\pi_{1}(U)\longrightarrow\pi_{1}(R)$ be the induced 
homomorphism of fundamental groups. We set
$K(U):=Ker(i_{*})\subset\pi_{1}(U)$.
Consider the regular covering $p_{U}:\tilde U\longrightarrow U$ of $U$ 
corresponding to the group $K(U)$, that is, $\pi_{1}(\tilde U)=K(U)$ and
$\pi_{1}(U)/K(U)$ acts on $\tilde U$ as the group of deck transformations.
Further, by 
$p_{U}^{*}(H^{\infty}(U))\subset H^{\infty}(\tilde U)$ we denote the subspace
of holomorphic functions invariant with respect to the action of 
$\pi_{1}(U)/K(U)$ (i.e. the pullback by $p_{U}$ of 
$H^{\infty}(U)$ to $\tilde U$).
Let $F_{z}:=p_{U}^{-1}(z)$, $z\in U$, and $l^{\infty}(F_{z})$
be the Banach space of bounded complex-valued functions on $F_{z}$ with
the supremum norm. By $C(F_{z})\subset l^{\infty}(F_{z})$ we denote
the subspace of constant functions.
\begin{Th}\label{fortype1}
There is a linear continuous projector
$P:H^{\infty}(\tilde U)\longrightarrow p_{U}^{*}(H^{\infty}(U))$ satisfying
the properties:
\begin{description}
\item[{\rm (1)}] \ There exists a family of linear continuous
projectors $P_{z}:l^{\infty}(F_{z})\longrightarrow C(F_{z})$
holomorphically depending on $z\in U$ such that 
$P[f]|_{p_{U}^{-1}(z)}:=P_{z}[f|_{p_{U}^{-1}(z)}]$ for any 
$f\in H^{\infty}(\tilde U)$;
\item[{\rm (2)}] \ $P(fg)=P(f)g$ for any $f\in H^{\infty}(\tilde U),\
g\in p_{U}^{*}(H^{\infty}(U))$;
\item[{\rm (3)}] \ If $f\in H^{\infty}(\tilde U)$ is such that 
$f|_{F_{z}}$ is constant, then $P(f)|_{F_{z}}=f|_{F_{z}}$;
\item[{\rm (4)}] \ Each $P_{z}$ is continuous in the weak $*$ topology of 
$l^{\infty}(F_{z})$;
\item[{\rm (5)}] \ The norm $||P||\leq C<\infty$ where $C=C(N)$ 
depends on $N$ only.
\end{description}
\end{Th}
As a simple corollary of Theorem \ref{fortype1} we obtain
\begin{C}\label{cor1}
Let $\tilde R$ be a covering of a
bordered Riemann surface $R$. The fundamental group $\pi_{1}(\tilde R)$ is a 
free group
whose family of generators $J$ is finite or countable. Assume
that $U\subset\tilde R$ is a domain such that $\pi_{1}(U)$
is generated by a subfamily of $J$.
Let $r:\Di\longrightarrow U$ 
be the universal covering map. Then there exists a linear continuous 
projector 
$P:H^{\infty}(\Di)\longrightarrow r^{*}(H^{\infty}(U))$ satisfying
the properties of Theorem \ref{fortype1}.
\end{C}
\begin{E}\label{e1}
{\rm Let $r:\Di\longrightarrow X$ be the universal covering of a
compact complex Riemann surface of genus $g\geq 2$. Let $K\subset\Di$ be the 
{\em fundamental compact} with respect to the action of the
deck transformation group $\pi_{1}(X)$. By definition, the boundary of
$K$ is the union of $2g$ analytic curves. Let 
$D_{1},...,D_{k}$ be a family of mutually disjoint closed disks situated
in the interior of $K$. We set 
$$
S:=\bigcup_{i=1}^{k}D_{i},\ \ \ K':=K\setminus S,\ \ \ {\rm and}\ \ \
\tilde R:=\bigcup_{g\in\pi_{1}(X)}g(K') .
$$
Then $R:=r(K')\subset X$ is a bordered Riemann surface, and 
$r:\tilde R\longrightarrow R$ is a regular covering corresponding to the
quotient group $\pi_{1}(X)$ of $\pi_{1}(R)$. Here 
$\pi_{1}(\tilde R)$ is generated by a family of simple closed curves in 
$\tilde R$ with the origin at
a fixed point $x_{0}\in\tilde R$ so that each such curve 
goes around only of one of $g(D_{i})$, 
$g\in\pi_{1}(X)$, $i=1,...,k$. Let $Y\subset\Di$ be 
a simply connected domain with the property: there is a subset 
$L\subset\pi_{1}(X)$ so that 
$$
Y\bigcap\left(\bigcup_{g\in\pi_{1}(X)}g(S)\right)=\bigcup_{g\in L}g(S).
$$ 
Clearly $U:=Y\setminus(\cup_{g\in L}\ g(S))$ satisfies 
the conditions of Corollary \ref{cor1}. Therefore the projector $P$,
described above, exists for $U$.}
\end{E}
\begin{R}\label{conj}
{\rm In view of Example \ref{e1} it is natural to conjecture the following.
\\
Let $U\subset\Di$ be a domain obtained by removing from $\Di$ a finite or 
countable family of pairwise disjoint closed disks $D_{i}$. 
Let $r_{i}$ be the radius of $D_{i}$ with respect to the pseudohyperbolic 
metric $\rho$,} 
\begin{equation}\label{pseudometric}
\rho(z,w):=\left|\frac{z-w}{1-\overline{z}w}\right|,\ \ \ z,w\in\Di\ .
\end{equation}
{\rm Assume that $\inf_{i}r_{i}=c>0$. Let 
$r:\Di\longrightarrow U$ be the universal covering map.\\
{\bf Conjecture.} {\em There is a linear continuous projector 
$P:H^{\infty}(\Di)\longrightarrow r^{*}(H^{\infty}(U))$ satisfying
$P(fg)=P(f)g$ for $f\in H^{\infty}(\Di)$, $g\in r^{*}(H^{\infty}(U))$ whose
norm depends on $c$ only.}}
\end{R}
To formulate our next corollary we recall several definitions.

Let $X$ be a Riemann surface such that $H^{\infty}(X)$ separates points of
$X$. By 
$M(H^{\infty}(X))$ we denote the maximal ideal space of $H^{\infty}(X)$, i.e. 
the set of nontrivial
multiplicative linear functionals on $H^{\infty}(X)$ with the weak $*$ 
topology (which is called the {\em Gelfand topology}).
It is a compact Hausdorff space. Each point $x\in X$ corresponds in a
natural way (point evaluation) to an element of $M(H^{\infty}(X))$.
So $X$ is naturally embedded into $M(H^{\infty}(X))$. Then the corona
problem for $H^{\infty}(X)$ asks: Is $M(H^{\infty}(X))$ the closure (in
the Gelfand topology) of $X$? 

Recall also that the corona problem has the following analytic reformulation.

A collection $f_{1},...,f_{n}$ of functions from 
$H^{\infty}(X)$ satisfies the corona condition if
\begin{equation}\label{coro1}
|f_{1}(x)|+|f_{2}(x)|+...+|f_{n}(x)|\geq\delta>0\ \ \ {\rm for\ all}\
x\in X.
\end{equation}
The corona problem being solvable means that the Bezout equation
$$
f_{1}g_{1}+f_{2}g_{2}+...+f_{n}g_{n}\equiv 1
$$
has a solution $g_{1},...,g_{n}\in H^{\infty}(X)$ for any 
$f_{1},...,f_{n}$ satisfying the corona condition. We refer to
$\max_{j}||g_{j}||$ as a ``bound on the corona solutions''.
Using Carleson's solution [Ca] of the corona problem for 
$H^{\infty}(\Di)$ and property
(2) for the projector constructed in Theorem \ref{fortype1}
we immediately obtain.
\begin{C}\label{cor2}
Let $N\subset\subset M$, $R\in {\cal F}_{c}(N)$ and $i:U\hookrightarrow R$
be open Riemann surfaces satisfying the conditions of
Theorem \ref{fortype1}. Assume also that $K(U):=Ker(i_{*})$ is 
trivial. Let $f_{1},...,f_{n}\in H^{\infty}(U)$
satisfy the corona condition (\ref{coro1}). Then the corona problem
has a solution $g_{1},...,g_{n}\in H^{\infty}(U)$ with the bound
$\max_{j}||g_{j}||\leq C(N,n,\delta/\max_{j}||f_{j}||)$.
\end{C}
\begin{R}\label{re1}
{\rm In this case $\pi_{1}(N)$ and $\pi_{1}(M)$ are free groups. Therefore 
condition (\ref{condit1}) is equivalent to $H_{1}(M,\Z)=H_{1}(N,\Z)$
for the corresponding homology groups.} 
\end{R}
{\bf 1.2.} Another application of Theorem \ref{fortype1} is a result on
the classification of interpolating sequences in $U$ (cf. [St] and [JM]).
Recall that a sequence $\{z_{j}\}\subset U$ 
is an {\em interpolating} sequence for $H^{\infty}(U)$ if
for every bounded sequence of complex numbers $\{a_{j}\}$, there is an
$f\in H^{\infty}(U)$ so that $f(z_{j})=a_{j}$.
\begin{Th}\label{inter}
Let $N\subset\subset M$, $R\in {\cal F}_{c}(N)$, $i:U\hookrightarrow R$ and 
$\tilde U$ be
complex manifolds satisfying the 
conditions of Theorem \ref{fortype1}. A sequence
$\{z_{j}\}\subset U$ is interpolating for $H^{\infty}(U)$ if and only if 
$r^{-1}(\{z_{j}\})$ is interpolating for $H^{\infty}(\tilde U)$.
\end{Th}
\begin{E}\label{e2}
{\rm Let $M\subset\Di$ be a bounded domain, whose boundary $B$ consists 
of $k$ simple closed continuous curves $B_{1},...,B_{k}$, with
$B_{1}$ forming the outer boundary. Let $D_{1}$ be the interior of $B_{1}$,
and $D_{2},...,D_{k}$ the exteriors of $B_{2},...,B_{k}$, including the point
at infinity. Then each $D_{i}$ is biholomorphic to $\Di$. Let
$\{z_{ji}\}_{j\in J}$ be an interpolating sequence for $H^{\infty}(D_{i})$, 
$i=1,...,k$,
such that
the Euclidean distance between any two distinct sequences is bounded
from below by a positive number. Then for any covering $p:R\longrightarrow M$
the sequence $p^{-1}(\{z_{ji}\}_{i,j})$ is interpolating for $H^{\infty}(R)$.}
\end{E}
In Section 5 we also establish some results for interpolating
sequences in $U$ with $U$ being a Riemann surface satisfying the assumptions
of Corollary \ref{cor2}. These results have much in common with similar
properties of interpolating sequences for $H^{\infty}(\Di)$.
%==============================================
\sect{\hspace*{-1em}. Construction of Bundles.}
In this section we formulate and prove some preliminary results used in the 
proofs of our main theorems. \\
{\bf 2.1. Definitions and Examples.} 
(For standard facts about bundles see e.g. Hirzebruch's book [Hi].) 
In what follows all topological spaces are assumed to be finite or 
infinite dimensional.

Let $X$ be a complex analytic space and $S$ be a complex 
analytic Lie group with
the unit $e\in S$. Consider an effective holomorphic action of $S$ on
a complex analytic space $F$. Here {\em holomorphic action} means a 
holomorphic map $S\times F\longrightarrow F$ sending 
$s\times f\in S\times F$ to
$sf\in F$ such that $s_{1}(s_{2}f)=(s_{1}s_{2})f$ and $ef=f$ for any
$f\in F$. {\em Efficiency} means that the condition $sf=f$ for some $s$ and 
any $f$ implies that $s=e$. 
\begin{D}\label{de1}
A complex analytic space $W$ together with a holomorphic map (projection)
$\pi:W\longrightarrow X$ is called a holomorphic bundle over $X$ with the
structure group $S$ and the fibre $F$, if there exists a
system of coordinate transformations, i.e., if
\\
(1) there is an open cover ${\cal U}=\{U_{i}\}_{i\in I}$ of $X$ and a
family of biholomorphisms 
$h_{i}:\pi^{-1}(U_{i})\longrightarrow U_{i}\times F$,
that map ``fibres'' $\pi^{-1}(u)$ onto $u\times F$;
\\
(2) for any $i,j\in I$ there are elements 
$s_{ij}\in {\cal O}(U_{i}\cap U_{j}, S)$ such that
$$
(h_{i}h_{j}^{-1})(u\times f)=u\times s_{ij}(u)f\ \ \ {\rm for\ any}\ \
u\in U_{i}\cap U_{j},\ f\in F\ .
$$
In particular, a 
holomorphic bundle $\pi:W\longrightarrow X$ whose fibre is a Banach space 
$F$ and the structure group is
$GL(F)$ (the group of linear invertible transformations of $F$) is
called a holomorphic Banach vector bundle. 

A holomorphic section of a holomorphic bundle $\pi:W\longrightarrow X$
is a holomorphic map $s:X\longrightarrow W$ satisfying $\pi\circ s=id$.
Let $\pi_{i}:W_{i}\longrightarrow X$, $i=1,2$, be holomorphic Banach vector
bundles. A holomorphic map $f:W_{1}\longrightarrow W_{2}$ satisfying\\
(a) $f(\pi_{1}^{-1}(x))\subset\pi_{2}^{-1}(x)$ for any $x\in X$;\\
(b) $f|_{\pi_{1}^{-1}(x)}$ is a linear continuous map of the corresponding
Banach spaces,\\
is called a homomorphism. If, in addition, $f$ is a 
homeomorphism, then $f$ is called an isomorphism.
\end{D}
We also use
the following construction of holomorphic bundles 
(see, e.g. [Hi,Ch.1]):

Let $S$ be a complex analytic Lie group and 
${\cal U}=\{U_{i}\}_{i\in I}$ be an open cover of $X$. 
By $Z_{\cal O}^{1}({\cal U},S)$ we denote the set of holomorphic $S$-valued 
${\cal U}$-cocycles. By definition, 
$s=\{s_{ij}\}\in Z_{\cal O}^{1}({\cal U},S)$, 
where $s_{ij}\in {\cal O}(U_{i}\cap U_{j}, S)$ and
$s_{ij}s_{jk}=s_{ik}\ \ \  {\rm on}\ \ \  U_{i}\cap U_{j}\cap U_{k}$.
Consider disjoint union $\sqcup_{i\in I}U_{i}\times F$ and for any 
$u\in U_{i}\cap U_{j}$ identify point $u\times f\in U_{j}\times F$ with 
$u\times s_{ij}(u)f\in U_{i}\times F$.
We obtain a holomorphic bundle $W_{s}$ over $X$ whose projection is induced 
by the projection $U_{i}\times F\longrightarrow U_{i}$.
Moreover, any holomorphic bundle over $X$ with the structure group $S$ and 
the fibre $F$ is isomorphic (in the category of holomorphic bundles) to a 
bundle $W_{s}$. 
\begin{E}\label{constbun}
{\rm {\bf (a)}  Let $M$ be a complex manifold.
For any subgroup $G\subset\pi_{1}(M)$ consider the unbranched covering
$g:M_{G}\longrightarrow M$ corresponding to $G$. We will describe $M_{G}$
as a holomorphic bundle over $M$.

First, assume that $G\subset\pi_{1}(M)$ is a normal subgroup.
Then $M_{G}$ is a regular covering of $M$ and the quotient group
$Q:=\pi_{1}(M)/G$ acts holomorphically on $M_{G}$ by deck 
transformations.  It is well known that $M_{G}$ in this case 
can be thought of as a principle fibre bundle over $M$ with fibre
$Q$ (here $Q$ is equipped with discrete topology). 
Namely, let us consider the map $R_{Q}(g):Q\longrightarrow Q$ 
defined by the formula
$$
R_{Q}(g)(h)=h\cdot g^{-1},\ \ \ h\in Q.
$$ 
Then there is an open cover ${\cal U}=\{U_{i}\}_{i\in I}$ of $M$ by sets
biholomorphic to open Euclidean balls in some $\Co^{n}$ and a locally constant
cocycle $c=\{c_{ij}\}\in Z_{\cal O}^{1}({\cal U},Q)$ 
such that $M_{G}$ is biholomorphic
to the quotient space of the disjoint union 
$V=\sqcup_{i\in I}U_{i}\times Q$ by the equivalence relation:
$U_{i}\times Q\ni x\times R_{Q}(c_{ij})(h)\sim x\times h\in 
U_{j}\times Q$.
The identification space is a holomorphic bundle with projection 
$p:M_{G}\longrightarrow M$ induced by the
projections $U_{i}\times Q\longrightarrow U_{i}$.
In particular, when $G=e$ we obtain the definition of the universal
covering $M_{e}$ of $M$.

Assume now that $G\subset\pi_{1}(M)$ is not necessarily normal.
Let $X_{G}=\pi_{1}(M)/G$ be the set of cosets with respect
to the (left) action of $G$ on $\pi_{1}(M)$ defined by left multiplications. 
By $[Gq]\in X_{G}$ we denote the coset containing $q\in\pi_{1}(M)$.
Let $H(X_{G})$ be the group of all homeomorphisms of $X_{G}$
(equipped with discrete topology). We define the homomorphism
$\tau:\pi_{1}(M)\longrightarrow H(X_{G})$ by the formula:
$$
\tau(g)([Gq]):=[Gqg^{-1}],\ \ \ q\in\pi_{1}(M).
$$
Set $Q(G):=\pi_{1}(M)/Ker(\tau)$ and let 
$\tilde g$ be the image of $g\in\pi_{1}(M)$ in $Q(G)$. 
By $\tau_{Q(G)}:Q(G)\longrightarrow H(X_{G})$ we denote the unique 
homomorphism whose pullback to $\pi_{1}(M)$ coincides with $\tau$.
Consider the action of
$G$ on $V=\sqcup_{i\in I}U_{i}\times \pi_{1}(M)$ induced by the left 
action of $G$ on $\pi_{1}(M)$ and let
$V_{G}=\sqcup_{i\in I}U_{i}\times X_{G}$ be the corresponding quotient set.
Define the equivalence relation
$U_{i}\times X_{G}\ni x\times \tau_{Q(G)}(\tilde c_{ij})(h)\sim x\times h\in 
U_{j}\times X_{G}$ 
with the same $\{c_{ij}\}$ as in the definition of $M_{e}$.
The corresponding quotient space is a holomorphic bundle with fibre $X_{G}$ 
biholomorphic to $M_{G}$. \\
{\bf (b)}\ We retain the notation of example (a). Let 
$B$ be a complex Banach space with norm $|\cdot|$.
Let $Iso(B)\subset GL(B)$ be the group
of linear isometries of $B$. Consider a homomorphism 
$\rho: Q\longrightarrow Iso(B)$. Without loss of generality we assume that 
$Ker(\rho)=e$, for otherwise we can pass to the corresponding quotient
group. The {\em holomorphic Banach vector bundle
$E_{\rho}\longrightarrow M$ associated with $\rho$ } is defined as 
the quotient of $\sqcup_{i\in I} U_{i}\times B$ by the equivalence
relation
$U_{i}\times B\ni x\times\rho(c_{ij})(w)\sim x\times w\in U_{j}\times B$
for any $x\in U_{i}\cap U_{j}$. Further, 
we can define a function $E_{\rho}\longrightarrow\Re_{+}$ which
will be called the {\em norm} on $E_{\rho}$ (and denoted by
the same symbol $|\cdot|$). The construction is as follows.  For any
$x\times w\in U_{i}\times B$ we set $|x\times w|:=|w|$. Since the image
of $\rho$ belongs to $Iso(B)$, the above definition is invariant with respect
to the equivalence relation determining $E_{\rho}$ and so it determines
a ``norm'' on $E_{\rho}$. Let us consider some examples. 

Let $l_{1}(Q)$ be the Banach space of complex-valued
sequences on $Q$ with $l_{1}$-norm. The action
$R_{Q}$ from (a) induces the homomorphism 
$\rho: Q\longrightarrow Iso(l_{1}(Q))$, 
$$
\rho(g)(w)[x]:=w(R_{Q}(g)(x)),\ \ \ g,x\in Q,\ w\in l_{1}(Q)\ .
$$
By $E_{1}^{M}(Q)$ we denote the holomorphic Banach vector bundle 
associated with $\rho$.

Let $l_{\infty}(Q)$ be the Banach space of bounded 
complex-valued sequences on $Q$ with $l_{\infty}$-norm.
The homomorphism $\rho^{*}: Q\longrightarrow Iso(l_{\infty}(Q))$, dual to 
$\rho$ is defined as 
$$
\rho^{*}(g)(v)[x]:=v(x\cdot g^{-1}),\ \ \ g,x\in Q,\ v\in l_{\infty}(Q).
$$ 
(It coincides with the homomorphism $(\rho^{t})^{-1}$:
$((\rho^{t})^{-1}(g)[v])(w):=v(\rho(g^{-1})[w])$,\
$g\in Q,\ v\in l_{\infty}(Q)$,\ $w\in l_{1}(Q)$.)
The holomorphic Banach vector bundle associated with $\rho^{*}$ will
be denoted by $E_{\infty}^{M}(Q)$. By definition it is dual to 
$E_{1}^{M}(Q)$.}
\end{E}
{\bf 2.2. Main Construction.} 
Let $B$ be a complex Banach space with norm $|\cdot |$ and let
$||\cdot||$ denote the corresponding norm on $GL(B)$.
For a discrete set $X$, denote by $B_{\infty}(X)$ the Banach space of
``sequences'' $b:=\{(x,b(x))\}_{x\in X}$, $b(x)\in B$, with norm
$$
|b|_{\infty}:=\sup_{x\in X}|b(x)| .
$$
By definition for 
$b_{i}=\{(x,b_{i}(x))\}_{x\in X}$, $\alpha_{i}\in\Co$, $i=1,2$, we have
$$
\alpha_{1}b_{1}+\alpha_{2}b_{2}=
\{(x,\alpha_{1}b_{1}(x)+\alpha_{2}b_{2}(x))\}_{x\in X} .
$$
Further, recall that a $B$-valued function
$f:U\longrightarrow B$ defined in an open set $U\subset\Co^{n}$
is said to be holomorphic if $f$ satisfies the
$B$-valued Cauchy integral formula in any polydisk containing in $U$.
Equivalently, locally $f$ can be represented as sum of absolutely
convergent holomorphic power series with coefficients in $B$.
Now any family $\{(x,f_{x})\}_{x\in X}$, where $f_{x}$ is a $B$-valued 
holomorphic on $U$ function
satisfying $|f_{x}(z)|<A$ for any $z\in U$ and
$x\in X$, can be considered as a $B_{\infty}(X)$-valued holomorphic function
on $U$.
In fact, the local Taylor expansion in this case follows from 
the Cauchy estimates of the coefficients in the 
Taylor expansion of each $f_{x}$.

Let  $t:X\longrightarrow X$ be a bijection and 
$h:X\times X\longrightarrow gl(B)$ be such that 
$$
h(t(x),x)\in GL(B)\ \ \ {\rm and}\ \ \
\max\{\sup_{x\in X}||h(t(x),x)||,\sup_{x\in X}||h^{-1}(x,t(x))||\}<\infty .
$$ 
Then we can define $a(h,t)\in GL(B_{\infty}(X))$ by the formula
$$
a(h,t)[(x,b(x))]:=(t(x),h(t(x),x)[b(x)]),\ \ \ 
b=\{(x,b(x))\}_{x\in X}\in B_{\infty}(X).
$$

We retain the notation of Example \ref{constbun}. For the acyclic cover
${\cal U}=\{U_{i}\}_{i\in I}$ of $M$ we have
$g^{-1}(U_{i})=\sqcup_{s\in X_{G}}V_{is}\subset M_{G}$ where 
$g|_{V_{is}}:V_{is}\longrightarrow U_{i}$ is biholomorphic.
Consider a holomorphic Banach vector bundle $\pi:E\longrightarrow M_{G}$
with fibre $B$ defined by 
coordinate transformations subordinate to the cover 
$\{V_{is}\}_{i\in I,s\in X_{G}}$ of $M_{G}$, i.e. by a holomorphic
cocycle $h=\{h_{is,jk}\}\in Z_{\cal O}^{1}(g^{-1}({\cal U}),GL(B))$,
$h_{is,jk}\in {\cal O}(V_{is}\cap V_{jk},GL(B))$, such that
$E$ is biholomorphic to the quotient space of disjoint union
$\sqcup_{i,s}V_{is}\times B$ by the equivalence relation
$V_{is}\times B\ni x\times h_{is,jk}(x)[v]\sim x\times v\in V_{jk}\times B$,
$s:=\tau_{Q(G)}(\tilde c_{ij})(k)$. The projection $\pi$ is
induced by coordinate projections $V_{is}\times B\longrightarrow V_{is}$.
Assume also that for any $x$
\begin{equation}\label{uniform}
\sup_{i,j,s,k}\max\{||h_{is,jk}(x)||,||h_{is,jk}^{-1}(x)||\}\leq A<\infty .
\end{equation}
Further, define $\tilde\pi:=g\circ\pi :E\longrightarrow M$.
\begin{Proposition}\label{pr1}
The triple $(E, M,\tilde\pi)$ determines a holomorphic Banach
vector bundle over $M$ with fibre $B_{\infty}(X_{G})$. (We denote this
bundle by $E_{M}$.)
\end{Proposition}
{\bf Proof.} 
Let $\phi_{is}:U_{i}\longrightarrow V_{is}$ be the map inverse to 
$g|_{V_{is}}$. We identify $V_{is}\times B$ with 
$U_{i}\times s\times B$ by $\phi_{is}$, and
$\{s\times B\}_{s\in X_{G}}$ with $B_{\infty}(X_{G})$. 
Further, for any $x\in U_{i}\cap U_{j}$, 
we set $\tilde h_{is,jk}(x):=h_{is,jk}(\phi_{is}(x))$.
Then
$E$ can be defined as the quotient space of
$\sqcup_{i\in I}U_{i}\times B_{\infty}(X_{G})$ by the equivalence 
relation 
$U_{j}\times B_{\infty}(X_{G})\ni x\times\{(k,b(k))\}_{k\in X_{G}}\sim
\{(\tau_{Q(G)}(\tilde c_{ij})(k),
\tilde h_{i\tau_{Q(G)}(\tilde c_{ij})(k),jk}(x)[b(k)])\}_{k\in X_{G}} 
\in U_{i}\times B_{\infty}(X_{G})$.
\\
Define
$\tilde h_{ij}(x):X_{G}\times X_{G}\longrightarrow gl(B)$, 
$x\in U_{i}\cap U_{j}$,
and $d_{ij}\in {\cal O}(U_{i}\cap U_{j},GL(B_{\infty}(X_{G})))$ by
the formulas
$$
\begin{array}{c}
\tilde h_{ij}(x)(s,k):=\tilde h_{is,jk}(x),\ \ \ {\rm and}\\
\\
d_{ij}(x)[b]:=a(\tilde h_{ij}(x),\tau_{Q(G)}(\tilde c_{ij}))[b],\ \ \ 
b\in B_{\infty}(X_{G}) .
\end{array}
$$
Here holomorphy of $d_{ij}$ follows from (\ref{uniform}). Clearly,
$d=\{d_{ij}\}$ is a holomorphic cocycle with values in 
$GL(B_{\infty}(X_{G}))$,
because $\{h_{is,jk}\}$ and $\{\tau_{Q(G)}(\tilde c_{ij})\}$ are cocycles.
Now $E$ can be
considered as a holomorphic Banach vector bundle over $M$ with fibre
$B_{\infty}(X_{G})$ obtained by identification in
$\sqcup_{i\in I}U_{i}\times B_{\infty}(X_{G})$ of 
$x\times d_{ij}(x)[b]\in U_{i}\times B_{\infty}(X_{G})$ with
$x\times b\in U_{j}\times B_{\infty}(X_{G})$, $x\in U_{i}\cap U_{j}$.
Moreover, according to our construction the projection 
$E\longrightarrow M$ coincides
with $\tilde\pi$.\ \ \ \ \ $\Box$

Let $W$ be a holomorphic Banach vector bundle
over a complex analytic space $X$. In what follows by
${\cal O}(U,W)$ we denote the vector space of holomorphic sections of
$W$ defined in an open set $U\subset X$.

We retain the notation of Proposition \ref{pr1}.
According to the construction of Proposition \ref{pr1},
a fibre $\tilde\pi^{-1}(z)$, $z\in U_{i}$, of $E_{M}$ can be 
identified with $\prod_{s\in X_{G}}\pi^{-1}(\phi_{is}(z))$ such that
if also $z\in U_{j}$ then 
\begin{equation}\label{prod1}
\prod_{s\in X_{G}}\pi^{-1}(\phi_{is}(z))=
\prod_{s\in X_{G}}\pi^{-1}(\phi_{j\tau_{Q(G)}(\tilde c_{ji})(s)}(z))\ .
\end{equation}
We recall the following definitions.

Let $J_{q}$ be the set of sequences $(i_{0}s_{0},...,i_{q}s_{q})$
with $i_{t}\in I$, $s_{t}\in X_{G}$ for $t=0,...,q$. A family 
$$
f=\{f_{i_{0}s_{0},...,i_{q}s_{q}}\}_{(i_{0}s_{0},...,i_{q}s_{q})
\in J_{q}},\ \ \ f_{i_{0}s_{0},...,i_{q}s_{q}}\in 
{\cal O}(V_{i_{0}s_{0}}\cap ...\cap V_{i_{q}s_{q}},E),
$$ 
is called a $q$-cochain on the cover 
$g^{-1}({\cal U}):=\cup_{i\in I,s\in X_{G}}V_{is}$ of $M_{G}$ with 
coefficients in the sheaf of germs of holomorphic sections of $E$.
These cochains generate a complex vector space $C^{q}(g^{-1}({\cal U}),E)$.
In the trivialization which identifies $\pi^{-1}(V_{i_{0}s_{0}})$ with
$V_{i_{0}s_{0}}\times B$ any $f_{i_{0}s_{0},...,i_{q}s_{q}}$ is represented
by $b_{i_{0}s_{0},...,i_{q}s_{q}}\in 
{\cal O}(V_{i_{0}s_{0}}\cap ...\cap V_{i_{q}s_{q}},B)$.
Assume that for any $(i_{0}s_{0},...,i_{q}s_{q})\in J_{q}$ and any compact
$K\subset U_{i_{0}}\cap ...\cap U_{i_{q}}$ there is a constant
$C=C(K)$ such that
\begin{equation}\label{bound}
\sup_{s_{0},...,s_{q},z\in K}
|(b_{i_{0}s_{0},...,i_{q}s_{q}}\circ\phi_{i_{0}s_{0}})(z)|<C
\end{equation}
The set of cochains $f$ satisfying (\ref{bound}) is a vector subspace
of $C^{q}(g^{-1}({\cal U}),E)$ which will be denoted by
$C_{b}^{q}(g^{-1}({\cal U}),E)$. Further, the formula
\begin{equation}\label{delt}
(\delta^{q}f)_{i_{0}s_{0},...,i_{q+1}s_{q+1}}=\sum_{k=0}^{q+1}
(-1)^{k}r_{W}^{W_{k}}(f_{i_{0}s_{0},...,
\widehat{i_{k}s_{k}},...,i_{q+1}s_{q+1}}),
\end{equation}
where $f\in C^{q}(g^{-1}({\cal U}),E)$, determines a homomorphism
$$
\delta^{q}:C^{q}(g^{-1}({\cal U}),E)\longrightarrow 
C^{q+1}(g^{-1}({\cal U}),E).
$$
Here $\widehat{}$ over a symbol means that this symbol must be omitted.
Besides, we set
$W=V_{i_{0}s_{0}}\cap ...\cap V_{i_{q+1}s_{q+1}}$,\ \
$W_{k}=V_{i_{0}s_{0}}\cap ...\cap\widehat{V}_{i_{k}s_{k}}\cap...
\cap V_{i_{q+1}s_{q+1}}$ and  $r_{W}^{W_{k}}$
is restriction map from $W$ to $W_{k}$. Also condition (\ref{uniform}) 
implies that $\delta^{q}$ maps $C_{b}^{q}(g^{-1}({\cal U}),E)$ into
$C_{b}^{q+1}(g^{-1}({\cal U}),E)$. We will denote 
$\delta^{q}|_{C_{b}^{q}(g^{-1}({\cal U}),E)}$
by $\delta_{b}^{q}$. As usual, $\delta^{q+1}\circ\delta^{q}=0$ and
$\delta_{b}^{q+1}\circ\delta_{b}^{q}=0$. Thus one can define the
cohomology groups on the cover $g^{-1}({\cal U})$ by
$$
H^{q}(g^{-1}({\cal U}),E):=Ker(\delta^{q})/Im(\delta^{q-1}),\ \ \
H_{b}^{q}(g^{-1}({\cal U}),E):=Ker(\delta_{b}^{q})/Im(\delta_{b}^{q-1})\ .
$$
In what follows the cohomology group $H^{q}({\cal U},E_{M})$ on the
cover ${\cal U}$ of $M$ with coefficients in the sheaf of germs of
holomorphic sections of $E_{M}$ is defined similarly to
$H^{q}(g^{-1}({\cal U}),E)$. Elements of $Ker(\delta^{q})$ and
$Ker(\delta_{b}^{q})$ will be called $q$-cocycles and of
$Im(\delta^{q-1})$ and $Im(\delta_{b}^{q-1})$  
$q$-coboundaries.
\begin{Proposition}\label{isom}
There is a linear isomorphism $\Phi^{q}:H_{b}^{q}(g^{-1}({\cal U}),E)
\rightarrow H^{q}({\cal U},E_{M})$.
\end{Proposition}
{\bf Proof.}
Let $f=\{f_{i_{0}s_{0},...,i_{q}s_{q}}\}\in C_{b}^{q}(g^{-1}({\cal U}),E)$. 
Let $b_{i_{0}s_{0},...,i_{q}s_{q}}\in 
{\cal O}(V_{i_{0}s_{0}}\cap ...\cap V_{i_{q}s_{q}},B)$ be 
the representation of $f_{i_{0}s_{0},...,i_{q}s_{q}}$ in the trivialization 
which identifies $\pi^{-1}(V_{i_{0}s_{0}})$ with $V_{i_{0}s_{0}}\times B$.
If $V_{i_{0}s_{0}}\cap ...\cap V_{i_{q}s_{q}}\neq\emptyset$ then
$s_{k}=\tau_{Q(G)(\tilde c_{i_{k}i_{0}})}(s_{0})$, $k=0,...,q$, and
$U_{i_{0}}\cap...\cap U_{i_{q}}\neq\emptyset$. For otherwise,
$b_{i_{0}s_{0},...,i_{q}s_{q}}=0$.
Thus for $s_{0},...,s_{q}$ satisfying the above identities we can define
$$
\tilde b_{i_{0},...,i_{q}}:=
\{b_{i_{0}s_{0},...,i_{q}s_{q}}\circ\phi_{i_{0}s_{0}}\}_{s_{0}\in X_{G}} .
$$
For $U_{i_{0}}\cap...\cap U_{i_{q}}=\emptyset$ we set
$\tilde b_{i_{0},...,i_{q}}=0$.
Further, according to (\ref{bound}), 
$\tilde b_{i_{0},...,i_{q}}\in 
{\cal O}(U_{i_{0}}\cap...\cap U_{i_{q}},B_{\infty}(X_{G}))$. This implies
that
$$
\tilde f_{i_{0},...,i_{q}}:=
\{f_{i_{0}s_{0},...,i_{q}s_{q}}\circ\phi_{i_{0}s_{0}}\}_{s_{0}\in X_{G}}
$$ 
defined similarly to $\tilde b_{i_{0},...,i_{q}}$ belongs to
${\cal O}(U_{i_{0}}\cap...\cap U_{i_{q}},E_{M})$, because
$\tilde b_{i_{0},...,i_{q}}$ is just another representation of
$\tilde f_{i_{0},...,i_{q}}$ in the trivialization which identifies
$\pi^{-1}(U_{i_{0}})$ with $U_{i_{0}}\times B_{\infty}(X_{G})$.
For $\tilde f=\{\tilde f_{i_{0},...,i_{q}}\}$ we set
$\tilde\Phi^{q}(f)=\tilde f$. Then, clearly,
$\tilde\Phi^{q}:C_{b}^{q}(g^{-1}({\cal U}),E)\longrightarrow 
C^{q}({\cal U},E_{M})$ is linear and injective. Now for a cochain
$\tilde f\in C^{q}({\cal U},E_{M})$ we can convert the  
construction for $\tilde\Phi^{q}$ to find a cochain
$f\in C_{b}^{q}(g^{-1}({\cal U}),E)$ such that $\tilde\Phi^{q}(f)=\tilde f$.
Thus $\tilde\Phi^{q}$ is an isomorphism. Moreover, a simple calculation
based on (\ref{prod1}) shows that
\begin{equation}\label{comute}
\delta^{q}\circ\tilde\Phi^{q}=\tilde\Phi^{q+1}\circ\delta_{b}^{q},
\end{equation}
where $\delta^{q}$ on the left means the operator for $E_{M}$ defined
similarly to (\ref{delt}). Hence $\tilde\Phi^{q}$ determines a 
linear isomorphism $\Phi^{q}:H_{b}^{q}(g^{-1}({\cal U}),E)
\longrightarrow H^{q}({\cal U},E_{M})$.
\ \ \ \ \ $\Box$

We complete this section by 
\begin{Proposition}\label{isobun}
Let $\rho:G\longrightarrow Iso(B)$ be a homomorphism and 
$E_{\rho}\longrightarrow M_{G}$ be the holomorphic Banach vector bundle
associated with $\rho$. Then $E_{\rho}$ satisfies conditions of 
Proposition \ref{pr1}.
\end{Proposition}
{\bf Proof.} Let $M_{e}\longrightarrow M_{G}$ be the universal covering
(recall that $G=\pi_{1}(M_{G})$).
Since the open cover 
$g^{-1}({\cal U})=\{V_{is}\}_{i\in I,s\in X_{G}}$ of $M_{G}$ is
acyclic, $M_{e}$ can be defined with respect to $g^{-1}({\cal U})$.
Namely, there is a cocycle $h=\{h_{is,jk}\}\in 
Z_{\cal O}^{1}(g^{-1}({\cal U}),G)$ such that $M_{e}$ is biholomorphic
to the quotient space of $\sqcup_{i,s}V_{is}\times G$ by the equivalence
relation $V_{is}\times G\ni x\times R_{G}(h_{is,jk})(f)\sim x\times f\in
V_{jk}\times G$, $s=\tau_{Q(G)}(\tilde c_{ij})(k)$; here
$R_{G}(q)(f):=f\cdot q^{-1}$, $f,q\in G$. Now $E_{\rho}$ is biholomorphic
to the quotient space of $\sqcup_{i,s}V_{is}\times B$ by the equivalence
relation $V_{is}\times B\ni x\times\rho(h_{is,jk})(v)\sim x\times v\in
V_{jk}\times B$.  Clearly, the
family $\{\rho(h_{is,jk})\}$ satisfies estimate (\ref{uniform}).
\ \ \ \ \ $\Box$
%========================================
\sect{\hspace*{-1em}. Proof of Theorem \ref{fortype1} and
Corollaries \ref{cor1}, \ref{cor2}.}
{\bf Proof of Theorem \ref{fortype1}.}
Let $N\subset\subset M$ be an open connected subset of a connected 
Stein manifold $M$ satisfying
(\ref{condit1}). Let $G\subset\pi_{1}(M)$ be a subgroup. As before, by 
$M_{G},N_{G}$ we denote the covering spaces of $M$ and $N$ corresponding
to $G$. Then by the covering homotopy theorem (see e.g.
[Hu,Ch.III,Sect.16]), there is
a holomorphic embedding $N_{G}\hookrightarrow M_{G}$. Without loss of 
generality we regard $N_{G}$ as an open subset of $M_{G}$. 
Denote also by $g_{MG}:M_{G}\longrightarrow M$,
$g_{NG}:N_{G}\longrightarrow N$ the corresponding projections such that
$g_{MG}|_{N_{G}}=g_{NG}$.
Let $i:U\hookrightarrow N_{G}$ be a holomorhic embedding of a complex
connected manifold $U$. 
\begin{Lm}\label{surj}
It suffices to prove the theorem under the assumption that
homomorphism
$i_{*}:\pi_{1}(U)\longrightarrow G\ (=\pi_{1}(N_{G}))$ is surjective.
\end{Lm}
{\bf Proof.}
Assume that $G':=Im(i_{*})$ is a proper subgroup of
$G$. By $t:N_{G'}\longrightarrow N_{G}$ we denote the covering of
$N_{G}$ corresponding to $G'\subset G$. By definition,
$g_{NG}\circ t=g_{NG'}:N_{G'}\longrightarrow N$ is the covering
of $N$ corresponding to $G'\subset\pi_{1}(N)$.
Further, by the covering homotopy theorem there is a holomorphic
embedding $i':U\hookrightarrow N_{G'}$ such that $t\circ i'=i$,\
$Ker(i'_{*})=Ker(i_{*})$, and $i'_{*}:\pi_{1}(U)\longrightarrow G'\
(=\pi_{1}(N_{G'}))$ is surjective. Clearly, it suffices to prove the
theorem for $i'(U)\subset N_{G'}$. \ \ \ \ \ $\Box$

In what follows we assume that $i_{*}$ is surjective.
By $p_{U}:\tilde U\longrightarrow U$ we denote the regular covering of $U$ 
corresponding to $K(U):=Ker(i_{*})$, where
$\pi_{1}(\tilde U)=K(U)$. 
Consider the holomorphic Banach vector bundle 
$E_{1}^{M_{G}}(G)\rightarrow M_{G}$ associated with homomorphism
$\rho_{G}:G\longrightarrow Iso(l_{1}(G))$,
$[\rho_{G}(g)(v)](x):=v(x\cdot g^{-1})$, $v\in l_{1}(G), x, 
g\in G$ (see Example \ref{constbun} (b)). Since $i_{*}$ is surjective,
$E_{1}^{M_{G}}(G)|_{U}=E_{1}^{U}(G)$.

Let $K_{G}\subset l_{1}(G)$ be the kernel of the linear
functional $l_{1}(G)\ni\{v_{g}\}_{g\in G}\mapsto
\sum_{g\in G}v_{g}$. Then $K_{G}$ is invariant with respect to any
$\rho_{G}(g)$, $g\in G$. In particular, $\rho_{G}$ determines a
homomorphism $h_{G}:G\longrightarrow Iso(K_{G})$, $h_{G}(g)=
\rho_{G}(g)|_{K_{G}}$. Here we consider $K_{G}$ with the norm induced
by the norm of $l_{1}(G)$. Let $F_{G}\longrightarrow M_{G}$ be the 
holomorphic Banach vector bundle associated with $h_{G}$. Clearly,
$F_{G}$ is a subbundle of $E_{1}^{M_{G}}(G)$.
Further, the quotient 
bundle $C_{G}:=E_{1}^{M_{G}}(G)/F_{G}\longrightarrow M_{G}$ is the trivial 
flat vector bundle of complex rank 1. Indeed, it is associated with 
the quotient homomorphism: 
$\tilde h_{G}:G\longrightarrow\Co^{*}$,\ $\tilde h_{G}(g)(v+K_{G}):=
\rho_{G}(g)(v)+K_{G}$, $g\in G,v\in l_{1}(G)$, where $w+K_{G}$ is the
image of $w\in l_{1}(G)$  in the factor space 
$l_{1}(G)/K_{G}=\Co$. This homomorphism is trivial because 
$\rho_{G}(g)(v)-v\in K_{G}$ by definition.
Thus we have the short exact sequence
\begin{equation}\label{exact}
0\longrightarrow F_{G}\longrightarrow E_{1}^{M_{G}}(G)
\stackrel{k_{G}}{\longrightarrow}
C_{G}\longrightarrow 0\ .
\end{equation}
Our goal is to construct a holomorphic section
$I_{G}:C_{G}\longrightarrow E_{1}^{M_{G}}(G)$ (linear on the fibres) such 
that $k_{G}\circ I_{G}=id$. 
Then we will obtain the bundle decomposition
$E_{1}^{M_{G}}(G)=I_{G}(C_{G})\oplus F_{G}$.

Let $\{t_{s}\}_{s\in G}$ be a
standard basis of unit vectors in $l_{1}(G)$, 
$t_{s}(g)=\delta_{sg}, \ s,g\in G$. Define 
$A:\Co\longrightarrow l_{1}(G)$
by $A(c)=ct_{e}$, where $e\in G$ is the unit. Then $A$ is a linear 
operator of norm $1$. Now let us recall the construction of 
$E_{1}^{M_{G}}(G)$ given in Proposition \ref{isobun}.

Let $M_{e}\longrightarrow M_{G}$ be the universal covering.
Consider an open cover
$g_{MG}^{-1}({\cal U})=\{V_{G,is}\}_{i\in I,s\in X_{G}}$ of $M_{G}$ where
${\cal U}:=\{U_{i}\}_{i\in I}$ is an open cover of $M$ by complex 
balls, and $\cup_{s\in X_{G}}V_{G,is}=g^{-1}(U_{i})$. Then
there is a cocycle $c_{G}=\{c_{G,is,jk}\}\in 
Z_{\cal O}^{1}(g_{MG}^{-1}({\cal U}),G)$ such that 
$E_{1}^{M_{G}}(G)$ is biholomorphic
to the quotient space of $\sqcup_{i,s}V_{G,is}\times l_{1}(G)$ by the 
equivalence relation 
$V_{G,is}\times l_{1}(G)\ni x\times\rho_{G}(c_{G,is,jk})(v)\sim x\times v\in
V_{G,jk}\times l_{1}(G)$. The construction of $F_{G}$ is similar, the 
only difference is that in the above formula we take $h_{G}$ instead of
$\rho_{G}$.
The above constructions restricted to $V_{G,is}$ determine 
isomorphisms of holomorphic Banach vector bundles:
$e_{G,is}:E_{1}^{M_{G}}(G)|_{V_{G,is}}\longrightarrow 
V_{G,is}\times l_{1}(G)$,
$f_{G,is}:F_{G}|_{V_{G,is}}\longrightarrow V_{G,is}\times K_{G}$,
$c_{G,is}:C_{G}|_{V_{G,is}}\longrightarrow V_{G,is}\times\Co$. Then we
define $A_{G,is}:C_{G}\longrightarrow E_{1}^{M_{G}}(G)$ on $V_{G,is}$ as 
$e_{G,is}^{-1}\circ A'\circ c_{G,is}$, where $A'(x\times c):=x\times A(c)$,
$x\in V_{G,is}$, $c\in\Co$.
Clearly, $k_{G}\circ A_{G,is}=id$ on $V_{G,is}$. Thus
$B_{G,is,jk}:=A_{G,is}-A_{G,jk}:
C_{G}|_{V_{G,is}\cap V_{G,jk}}\longrightarrow F_{G}|_{V_{G,is}\cap V_{G,jk}}$
is a homomorphism of bundles of norm $\leq 2$ on each fibre
(here norms on $F_{G}$,
$C_{G}$ and $E_{1}^{M_{G}}(G)$ are defined as in Example \ref{constbun} (b)). 
We also use the following identification 
$Hom(C_{G},F_{G})\cong F_{G}$ 
(the last isomorphism is because $C_{G}$ is trivial and 
$Hom(\Co,K_{G})\cong\Co^{*}\otimes K_{G}=K_{G}$).
Further, according to Proposition \ref{isobun}, the holomorphic Banach
vector bundle $Hom(C_{G},F_{G})$ associated with the homomorphism
$\tilde h_{G}\otimes h_{G}:G\longrightarrow Iso(Hom(\Co,K_{G}))$ satisfies
conditions of Proposition \ref{pr1}. Therefore, by definition,
$B_{G}=\{B_{G,is,jk}\}$ is a holomorphic 1-cocylce with respect to
$\delta_{b}^{1}$ defined on the cover $g_{MG}^{-1}({\cal U})$.
By $\phi_{G,is}:U_{i}\longrightarrow V_{G,is}$ we denote the map inverse to 
$g_{MG}|_{V_{G,is}}$. Then we will prove
\begin{Lm}\label{mainestim}
There is 
$\tilde B_{G}=\{\tilde B_{G,is}\}\in C_{b}^{0}(g_{MG}^{-1}({\cal U}),F_{G})$,
$\tilde B_{G,is}\in {\cal O}(V_{is},F_{G})$,
such that $\delta_{b}^{0}(\tilde B_{G})=B_{G}$. Moreover, 
for any $i\in I$ there is a continuous nonnegative function 
$F_{i}:U_{i}\longrightarrow\Re_{+}$ such that for any $G$ 
\begin{equation}\label{afi}
\sup_{s\in X_{G},z\in U_{i}}|(\tilde B_{G,is}\circ\phi_{G,is})(z)|\leq
F_{i}(z)\ .
\end{equation}
Here $|\cdot|$ denotes the norm on $F_{G}$.
\end{Lm}
{\bf Proof}.
According to Proposition \ref{pr1}, we can construct the holomorphic
Banach vector bundle $(F_{G})_{M}$. It is defined on the cover
${\cal U}$ of $M$ by a cocycle
$d_{G}:=\{d_{G,ij}\}\in 
Z_{\cal O}^{1}({\cal U},Iso((K_{G})_{\infty}(X_{G})))$, where
$d_{G,ij}\in {\cal O}(U_{i}\cap U_{j},Iso((K_{G})_{\infty}(X_{G})))$.
Let $\tilde\Phi_{G}^{q}:C_{b}^{q}(g_{MG}^{-1}({\cal U}),F_{G})
\longrightarrow C^{q}({\cal U},(F_{G})_{M})$ be the isomorphism
intorduced in the proof of Proposition \ref{isom}.
Then $\tilde\Phi_{G}^{1}(B_{G}):=b_{G}=\{b_{G,ij}\}$ is a holomorphic 
1-cocycle 
with respect to $\delta^{1}$ defined on ${\cal U}$. Here
$b_{G,ij}\in {\cal O}(U_{i}\cap U_{j}, (F_{G})_{M})$, and 
$$
\sup_{i,j\in I,z\in M}|b_{G,ij}(z)|_{(F_{G})_{M}}\leq 2,
$$
where $|\cdot|_{(F_{G})_{M}}$ denotes the norm on $(F_{G})_{M}$.

Let ${\cal G}$ be the set of all subgroups $G\subset\pi_{1}(M)$.
We define the Banach space $K=\oplus_{G\in {\cal G}}(K_{G})_{\infty}(X_{G})$
such that $x=\{x_{G}\}_{G\in {\cal G}}$ belongs to $K$ if 
$x_{G}\in (K_{G})_{\infty}(X_{G})$ and
$$
|x|:=\sup_{G\in {\cal G}}|x_{G}|_{(K_{G})_{\infty}(X_{G})}<\infty\ ,
$$
where $|\cdot|_{(K_{G})_{\infty}(X_{G})}$ is the norm on 
$(K_{G})_{\infty}(X_{G})$. Further, let us define 
$d:=\{d_{ij}\}\in Z_{\cal O}^{1}({\cal U}, Iso(K))$ as
$d:=\oplus_{G\in {\cal G}}d_{G}$.
Here
$$
d_{ij}:=\oplus_{G\in {\cal G}}d_{G,ij},\ \ \
[d_{ij}(z)](\{v_{G}\}_{G\in {\cal G}})
:=\{[d_{G,ij}(z)](v_{G})\}_{G\in {\cal G}},\ \ \ z\in U_{i}\cap U_{j}.
$$
Clearly $d_{ij}\in {\cal O}(U_{i}\cap U_{j},Iso(K))$. Now we define
the holomorphic Banach vector bundle $F$ over $M$ by the identification
$U_{i}\times K\ni x\times d_{ij}(x)[v]\sim x\times v\in U_{j}\times K$
for any $x\in U_{i}\cap U_{j}$. In fact, this bundle coincides with
$\oplus_{G\in {\cal G}}(F_{G})_{M}$. A vector $f$ of $F$ over $z\in M$
can be identify with a family $\{f_{G}\}_{G\in {\cal G}}$ so that
$f_{G}\in (F_{G})_{M}$ is a vector over $z$. Moreover, 
the norm $|f|_{F}:=\sup_{G\in {\cal G}}|f_{G}|_{(F_{G})_{M}}$ of $f$ 
is finite.
Now we can define a holomorphic 1-cocycle $b=\{b_{ij}\}$ of $F$  defined on 
the cover ${\cal U}$ as
$$
b:=\{b_{G}\}_{G\in {\cal G}},\ \ \ 
b_{ij}:=\{b_{G,ij}\}_{G\in {\cal G}}\in {\cal O}(U_{i}\cap U_{j},F)\ .
$$
Here holomorphy of $b_{ij}$ follows from the uniform estimate of
norms of $b_{G,ij}$.

Let us use the fact that $M$ is a Stein manifold. According to a 
theorem of Bungart [B, Sect.4] (i.e. a version of the
classical Cartan Theorem B for cohomology of sheaves of germs of
holomorphic sections of holomorphic Banach vector bundles) cocycle
$b$ represents 0 in the corresponding cohomology
group $H^{1}(M,F)$. Further, the cover $\{U_{i}\}_{i\in I}$ of $M$ 
consists of Stein manifolds (and so it is acyclic). Therefore
by the  classical Ler\'{e} theorem (on calculation of cohomology
groups on acyclic covers),
$$
H^{1}(M,F)=H^{1}({\cal U},F)\ .
$$
Thus $b$ represents 0 in $H^{1}({\cal U},F)$, that is, $b$ is a
coboundary. In particular, there are holomorphic sections
$b_{i}\in {\cal O}(U_{i},F)$  such that
$$ 
b_{i}(z)-b_{j}(z)=b_{ij}(z)\ \ \ {\rm for\ any}\ z\in U_{i}\cap U_{j}\ .
$$
We also set 
$$
F_{i}(z):=|b_{i}(z)|_{F}\ .
$$
Then $F_{i}$ is a continuous nonnegative function on $U_{i}$.
Further, by definition each $b_{i}$ can be represented as a family
$\{b_{G,i}\}_{G\in {\cal G}}$ where 
$b_{G,i}\in {\cal O}(U_{i},(F_{G})_{M})$. The family 
$\tilde b_{G}=\{b_{G,i}\}_{i\in I}$ belongs to 
$C^{0}({\cal U}, (F_{G})_{M})$. Using the isomorphism 
$\tilde\Phi_{G}^{0}$ from
Proposition \ref{isom} we obtain a cochain
$\tilde B_{G}:=[\tilde\Phi_{G}^{0}]^{-1}(\tilde b_{G})\in 
C_{b}^{0}(g_{MG}^{-1}({\cal U}), F_{G})$. Now if 
$\tilde B_{G}:=\{\tilde B_{G,is}\}$, 
$\tilde B_{G,is}\in {\cal O}(V_{is},F_{G})$, from identity
(\ref{comute}) it follows that
$$
\tilde B_{G,is}(z)-\tilde B_{G,jk}(z)=B_{G,is,jk}(z)\ \ \ {\rm for\ any}\ 
z\in V_{is}\cap V_{jk}.
$$
Finally, inequality (\ref{afi}) is the consequence of definitions of
$F_{i}$ and $\tilde\Phi_{G}^{0}$.

The lemma is proved.\ \ \ \ \ $\Box$

Let us consider now the family $\{A_{G,is}-B_{G,is}\}_{i,s}$. 
By definition, it determines a
holomorphic linear section  $I_{G}:C_{G}\longrightarrow 
E_{1}^{M_{G}}(G)$, $k_{G}\circ I_{G}=id$. Thus we have 
$E_{1}^{M_{G}}(G)=I_{G}(C_{G})\oplus F_{G}$. In the next result the norm
$||\cdot||$ of $I_{G}$ is defined with respect to the norms 
$|\cdot|_{C_{G}}$ and $|\cdot|_{E_{1}^{M_{G}}(G)}$.
\begin{Lm}\label{normest}
There is a constant $C=C(N)$ such that for any $G\in {\cal G}$
$$
\sup_{z\in N}||I_{G}(z)||\leq C\ .
$$
\end{Lm}
{\bf Proof.} Let ${\cal V}=\{V_{i}\}$ be a refinement of the 
cover ${\cal U}$ of $M$ such that each $V_{i}$ is relatively compact in 
some $U_{k(i)}$. Then from Lemma \ref{mainestim} it follows that
$$
\sup_{s\in X_{G},z\in V_{i}}|(\tilde B_{G,k(i)s}\circ\phi_{G,k(i)s})(z)|\leq
\sup_{z\in V_{i}}F_{k(i)}(z)=C_{i}<\infty\ .
$$
Now for any $z\in g_{MG}^{-1}(V_{i})$ we have
$$
||I_{G}(z)||\leq
\sup_{s\in X_{G}, y\in V_{i}}\{||(A_{G,k(i)s}\circ\phi_{G,k(i)s})(y)||+
||(\tilde B_{G,k(i)s}\circ\phi_{G,k(i)s})(y)||\}
\leq 1+C_{i}\ .
$$
Since $\overline{N}\subset M$ is a compact, we can find a finite number
of sets $V_{i_{1}},...,V_{i_{l}}$ which cover $\overline{N}$. Then
$$
\sup_{z\in N}||I_{G}(z)||\leq\max_{1\leq t\leq l}\{1+C_{i_{t}}\}:=C<\infty\ .
\ \ \ \ \ \Box
$$

Consider now the restriction of exact sequence (\ref{exact}) to $U$. Using 
the identification $E_{1}^{M_{G}}(G)|_{U}\cong E_{1}^{U}(G)$ we obtain
$$
0\longrightarrow (F_{G})|_{U}\longrightarrow E_{1}^{U}(G)\longrightarrow 
(C_{G})|_{U}\longrightarrow 0 .
$$
Similarly, we have the dual sequence obtained by taken the dual 
bundles in the above sequence
$$
0\longrightarrow [(C_{G})|_{U}]^{*}\longrightarrow E_{\infty}^{U}(G)
\longrightarrow [(F_{G})|_{U}]^{*}\longrightarrow 0 .
$$
Let $C(G)$ be the space of constant functions in $l_{\infty}(G)$. 
By definition, $[(C_{G})|_{U}]^{*}$ is a subbundle of 
$E_{\infty}^{U}(G)$ of complex
rank 1 with fibre $C(G)$ associated with the trivial homomorphism
$G\longrightarrow Iso(C(G))$. Let 
$P_{U}:=[(I_{G})|_{U}]^{*}: 
E_{\infty}^{U}(G)\longrightarrow [(C_{G})|_{U}]^{*}$
be the homomorphism of bundles dual to $(I_{G})|_{U}$.
Then for any $z\in U$, $P_{U}(z)$ projects the fibre of 
$E_{\infty}^{U}(G)$ over $z$ onto the fibre of $[(C_{G})|_{U}]^{*}$ 
over $z$. Moreover, we have
\begin{equation}\label{project}
\sup_{z\in U}||P_{U}(z)||\leq C,
\end{equation}
where $||\cdot||$ is the dual norm defined with respect to
$|\cdot|_{E_{\infty}^{U}(G)}$ and $|\cdot|_{[(C_{G})|_{U}]^{*}}$. The
operator $P_{U}$ induces also a linear map $P_{U}': 
{\cal O}(U,E_{\infty}^{U}(G))\longrightarrow 
{\cal O}(U,[(C_{G})|_{U}]^{*})$,
$$
[P_{U}'(f)](z):=[P_{U}(z)](f(z)),\ \ \ \ f\in {\cal O}(U,E_{\infty}^{U}(G))\ .
$$
Further, any $f\in H^{\infty}(\tilde U)$ can be considered in a
natural way as a bounded holomorphic  section of the trivial bundle 
$\tilde U\times\Co\longrightarrow\tilde U$. This bundle satisfies
assumptions of Proposition \ref{isobun} (for $U$ instead of $M$). 
Furthermore, it easy to see that in this case the bundle 
$(\tilde U\times\Co)_{U}$ defined in 
Proposition \ref{pr1} coincides with $E_{\infty}^{U}(G)$.
Let $\Phi_{U}^{0}:H_{b}^{0}(g^{-1}({\cal U}),\tilde U\times\Co)
\longrightarrow H^{0}({\cal U},E_{\infty}^{U}(G))$ be the isomorphism
of Proposition \ref{isom}. (This is just the direct image map 
with respect to $p_{U}:\tilde U\longrightarrow U$.)
We define the Banach subspace
$S_{\infty}(U)\subset H^{0}({\cal U},E_{\infty}^{U}(G))$ with
norm $|\cdot|_{U}$ by the formula
$$
f\in S_{\infty}(U)\ \ \ \leftrightarrow\ \ \
|f|_{U}:=\sup_{z\in U}|f(z)|_{E_{\infty}^{U}(G)}<\infty\ .
$$
Clearly $\Phi_{U}^{0}$ maps $H^{\infty}(\tilde U)$ 
isomorphically onto $S_{\infty}(U)$. Moreover, 
$s_{U}:=\Phi_{U}^{0}|_{H^{\infty}(\tilde U)}$ is a linear isometry of
Banach spaces. By definition, the space 
$s_{U}(p_{U}^{*}(H^{\infty}(U)))$ coincides with
${\cal O}(U,[(C_{U})|_{U}]^{*})\cap S_{\infty}(U)$. Then according to
the definition of $s_{U}$ and inequality
(\ref{project}) the linear operator
$P:=s_{U}^{-1}\circ P_{U}'\circ s_{U}$ maps 
$H^{\infty}(\tilde U)$ onto
$p_{U}^{*}(H^{\infty}(U))$. According to our construction $P$ is 
a bounded projector satisfying (1). Here the required projector
$P_{z}: l^{\infty}(F_{z})\longrightarrow C(F_{z})$ can be naturally identified
with $P_{U}(z)$. Let now $f\in H^{\infty}(\tilde U)$ and 
$g\in p_{U}^{*}(H^{\infty}(U))$. Then by definition we have
$$
P[f\cdot g]|_{p_{U}^{-1}(z)}=P_{z}[(f\cdot g)|_{p_{U}^{-1}(z)}]=
P_{z}[f|_{p_{U}^{-1}(z)}]\cdot g|_{p_{U}^{-1}(z)}=
(P[f]\cdot g)|_{p_{U}^{-1}(z)}\ .
$$
Here we used that $g|_{p_{U}^{-1}(z)}$ is a constant and $P_{z}$ is a
linear operator. This implies (2). Property (3) follows from the
fact that $P_{z}$ is a projector onto $C(F_{z})$. Further, (4) is a
consequence of the fact that $P_{U}(z)$ is dual to $(I_{G})(z)$ 
and so $P_{z}$ 
is continuous in the weak $*$ topology of $l_{\infty}(F_{z})$. Finally,
the norm of $P$ coincides with $\sup_{z\in U}||P_{U}(z)||$. Thus
$||P||\leq C$ for $C$ as in (\ref{project}). This completes the proof of
(5).

The theorem is proved.\ \ \ \ \ $\Box$\\
{\bf Proof of Corollary \ref{cor1}.} First, remark that any bordered
Riemann surface $N$ admits an embedding to a Riemann surface $M$ such that
the pair $N\subset\subset M$ satisfies condition (\ref{condit1}). 
Let $\tilde R$ be a covering of $N$
and $i:U\hookrightarrow\tilde R$ be such that 
$\pi_{1}(U)$ is generated by a subfamily of generators of the free group 
$\pi_{1}(\tilde R)$. Then the homomorphism 
$i_{*}:\pi_{1}(U)\longrightarrow\pi_{1}(\tilde R)$ is injective. 
In particular, $K(U):=Ker(i_{*})=\{e\}$ and
$p_{U}:\tilde U\longrightarrow U$ is the universal covering. Since
$\tilde U$ is biholomorphic to $\Di$, the existence of the projector 
$P:H^{\infty}(\Di)\longrightarrow p_{U}^{*}(H^{\infty}(U))$ follows from
Theorem \ref{fortype1}.\ \ \ \ \ $\Box$ \\
{\bf Proof of Corollary \ref{cor2}.}
Let $N\subset\subset M$, $R\subset {\cal F}_{c}(N)$ and
$i:U\hookrightarrow R$ be open Riemann surfaces satisfying conditions
of Theorem \ref{fortype1}. Assume also that $K(U):=Ker(i_{*})=\{e\}$.
Let $p_{U}:\Di\longrightarrow U$ be the universal covering map.
Then there is a projector
$P:H^{\infty}(\Di)\longrightarrow p_{U}^{*}(H^{\infty}(U))$ with
properties (1)-(5) of Theorem \ref{fortype1}. Let 
$f_{1},...,f_{n}\in H^{\infty}(U)$ satisfy the corona condition
(\ref{coro1}) with $\delta>0$. Without loss of generality we
will assume also that $\max_{i}||f_{i}||_{H^{\infty}(U)}\leq 1$.
For $1\leq i\leq n$ we set  $h_{i}:=p_{U}^{*}(f_{i})$. Then $h_{1},...,h_{n}
\in H^{\infty}(\Di)$ satisfy the corona condition in $\Di$ (with the same
$\delta$). Also $\max_{i}||h_{i}||_{H^{\infty}(\Di)}\leq 1$. Now 
according to the solution of Carleson's Corona Theorem [Ca], there is a 
constant
$C(n,\delta)$ and some $g_{1},...,g_{n}\in H^{\infty}(\Di)$ satisfying
$\max_{i}||g_{i}||_{H^{\infty}(\Di)}\leq C(n,\delta)$ such that 
$\sum_{i=1}^{n}g_{i}h_{i}\equiv 1$.
Let us define $d_{i}\in H^{\infty}(U)$ by the formula
$$
p_{U}^{*}(d_{i}):=P[g_{i}],\ \ \ 1\leq i\leq n\ .
$$
Then property (2) for $P$ implies that $\sum_{i=1}^{n}d_{i}f_{i}\equiv 1$.
Moreover, $\max_{i}||d_{i}||_{H^{\infty}(U)}\leq C(N)\cdot C(n,\delta)$ where
$C(N)$ is the constant from Lemma \ref{normest}.

The proof of the corollary is complete.\ \ \ \ \ $\Box$
\begin{R}\label{matrix}
{\rm In a forthcoming paper we present the following generalization of 
Corollary \ref{cor2}.}\\
{\bf Theorem.} 
Let an open Riemann surface $U$ satisfy the conditions of Corollary
\ref{cor2}. Let $b$ be an $n\times k$ matrix, $k<n$, with entries
in $H^{\infty}(U)$. Assume that the corona condition (\ref{coro1}) is valid 
for the family of minors of $b$ of order $k$. Then there is an 
$n\times n$ matrix $\tilde b$ with entries in $H^{\infty}(U)$ which 
extends $b$ such that $det(\tilde b)\equiv 1$ on $U$.

{\rm The proof of the theorem is based on Theorem \ref{fortype1} and
a Grauert type theorem for ``holomorphic'' vector bundles defined on maximal 
ideal spaces (which are not usual manifolds !) of certain Banach 
algebras.}
\end{R}
%====================
\sect{\hspace*{-1em}. Proof of Theorem \ref{inter}.}
Let $N\subset\subset M$ be a relatively compact domain of a connected
Stein manifold $M$ satisfying (\ref{condit1}). For a
subgroup $G\subset\pi_{1}(M)$ we denote by
$g_{NG}:N_{G}\longrightarrow N$ and $g_{MG}:M_{G}\longrightarrow M$
the covering spaces of $M$ and $N$ corresponding to the group $G$
with $N_{G}\subset M_{G}$. Further, assume that 
$i:U\hookrightarrow N_{G}$ is a
holomorphic embedding of a complex connected manifold $U$, 
$K(U):=Ker(i_{*})\subset\pi_{1}(U)$, and $p_{U}:\tilde U\longrightarrow U$
is the regular covering of $U$ corresponding to $K(U)$.
As before, without loss of generality we may assume that homomorphism
$i_{*}:\pi_{1}(U)\longrightarrow G\ (=\pi_{1}(N_{G}))$ is surjective
(see arguments of Lemma \ref{surj}). Thus the deck transformation group
of $\tilde U$ is $G$.
We begin the proof of the theorem with the following
\begin{Proposition}\label{intfibre}
For any $z\in U$, the sequence 
$p_{U}^{-1}(z):=\{w_{s}\}_{s\in G}\subset\tilde U$ is
interpolating with respect to $H^{\infty}(\tilde U)$. Moreover,
let 
$$
M(z)=\sup_{||a_{s}||_{l^{\infty}(G)}\leq 1}
\inf\{||g||_{H^{\infty}(\tilde U)} :\ g\in H^{\infty}(\tilde U),\
g(w_{s})=a_{s},\ j=1,2,...\}
$$
be the constant of interpolation for $p_{U}^{-1}(z)$. Then there
is a constant $C=C(N)$ such that
$$
\sup_{z\in U}M(z)\leq C\ .
$$
\end{Proposition}
{\bf Proof.}
Consider the homomorphism 
$\rho_{G}^{*}: G\longrightarrow Iso(l_{\infty}(G))$,
$$
[\rho_{G}^{*}(g)(w)](x):=w(x\cdot g^{-1}),\ \ \ w\in l_{\infty}(G), x, 
g\in G\ .
$$ 
Let $E_{\infty}^{M_{G}}(G)\longrightarrow M_{G}$ be the holomorphic
Banach vector bundle associated with $\rho_{G}^{*}$. Then 
$E_{\infty}^{M_{G}}(G)|_{U}=E_{\infty}^{U}(G)$ (see Example \ref{constbun}
(b)). According to Proposition \ref{isobun}, we can define the 
holomorphic Banach vector bundle 
$[E_{\infty}^{M_{G}}(G)]_{M}\longrightarrow M$ with the fibre
$[l^{\infty}(G)]_{\infty}(X_{G})$.
Let ${\cal G}$ be the set of all subgroup $G\subset\pi_{1}(M)$.
We define the Banach space 
$L=\oplus_{G\in {\cal G}}[l^{\infty}(G)]_{\infty}(X_{G})$ such that
$x=\{x_{G}\}_{G\in {\cal G}}$ belongs to $L$ if 
$x_{G}\in [l^{\infty}(G)]_{\infty}(X_{G})$ and
$$
|x|_{L}:=
\sup_{G\in {\cal G}}|x_{G}|_{[l^{\infty}(G)]_{\infty}(X_{G})}<\infty\ ,
$$
where $|\cdot|_{[l^{\infty}(G)]_{\infty}(X_{G})}$ is the norm on
$[l^{\infty}(G)]_{\infty}(X_{G})$. Then similarly to the construction of
Lemma \ref{mainestim}, we can define the holomorphic Banach vector bundle
$B$ over $M$ with the fibre $L$ by the formula
$$
B:=\oplus_{G\in {\cal G}}[E_{\infty}^{M_{G}}(G)]_{M}\ .
$$
Note that the structure group of $B$ is $Iso(L)$. Therefore the norm
$|\cdot|_{L}$ induces a norm $|\cdot|_{B}$ on $B$ (see Example 
\ref{constbun} (b)).
Let $O\subset\subset M$ be a relatively compact domain containing
$\overline{N}$. Denote by $H^{\infty}(O,B)$ the Banach space of bounded 
holomorphic sections from ${\cal O}(O,B)$, that is,
$$
f\in H^{\infty}(O,B)\ \ \ \leftrightarrow\ \ \ 
||f||:=\sup_{z\in O}|f(z)|_{B}<\infty\ .
$$
For any $z\in O$ consider the restriction operator 
$r(z):H^{\infty}(O,B)\longrightarrow L$,
$$
r(z)[f]:=f(z),\ \ \ f\in H^{\infty}(O,B)\ .
$$ 
Then $r(z)$ is
a continuous linear operator with the norm $||r(z)||\leq 1$. Moreover, by a
theorem of Bungart (see [B, Sect.4]), for any $v\in L$ there is a 
section $f\in {\cal O}(M,B)$ such that $f(z)=v$. Since $O$ is relatively
compact in $M$, the restriction $f|_{O}$ belongs to $H^{\infty}(O,B)$. This
shows that $r(z)$ is surjective. For any
$v\in L$ we set $K_{v}(z):=r(z)^{-1}(v)\subset H^{\infty}(O,B)$.
The constant
$$
h(z):=\sup_{|v|_{L}\leq 1}\inf_{t\in K_{v}(z)}||t||
$$
will be called {\em the constant of interpolation} for $r(z)$.
We will show that
\begin{Lm}\label{fibest}
$$
\sup_{z\in \overline{N}}h(z)\leq C<\infty
$$
where $C$ depends on $N$ only.
\end{Lm}
{\bf Proof.} In fact it suffices to cover $\overline{N}$
by a finite number of open balls and prove the required inequality for
$z$ varying in each of these balls. Moreover, since $\overline{N}$ is a 
compact, for any $w\in\overline{N}$ it suffices to find an open neighbourhood
$U_{w}\subset O$ of $w$ such that $\{h(z)\}_{z\in U_{w}}$ is bounded from 
above by an absolute constant.

Let $w\in\overline{N}$. Without
loss of generality we may identify a small open neighbourhood of $w$ in $O$
with the open unit ball $B_{c}(0,1)\subset\Co^{n}$, $n=dim\ O$, such 
that $w$ corresponds to 0 in this identification. 
It is easy to see that $r(z)$, $z\in B_{c}(0,1)$, is the family of 
linear continuous operators holomorphic in $z$. Let $R:=1/4h(w)$.
Since $h(w)\geq 1$, $B_{c}(0,1)$ contains $B_{c}(0,R)$. For 
a $y\in B_{c}(0,R)$ consider the one dimensional complex 
subspace $l_{y}$ of $\Co^{n}$ containing $y$. Without loss of generality
we may identify $l_{y}\cap B_{c}(0,1)$ with the open unit disk 
$\Di\subset\Co$. With this identification, let 
$r(z):=\sum_{i=0}^{\infty}r_{i}z^{i}$
be the Taylor expansion of $r(z)$ in $\Di$. Here 
$r_{i}:H^{\infty}(O,B)\longrightarrow L$ is a
linear operator with the norm $||r_{i}||\leq 1$. The last estimate follows 
from the Cauchy estimates for derivatives of holomorphic functions. We 
also have $r_{0}:=r(0)$ (recall that $w=0$). 
Let  $v\in L$, $|v|_{L}\leq 1$. For $z<R$ we will 
construct
$v(z)\in H^{\infty}(O,B)$ which depends holomorphically  on $z$, such that
$||v(z)||\leq 8h(w)$ and $r(z)[v(z)]=v$.\\
Let $v(z)=\sum_{i=0}^{\infty}v_{i}z^{i}$.
Then we have the formal decomposition
$$
v=r(z)[v(z)]=\sum_{i=0}^{\infty}z^{i}\cdot\sum_{j=0}^{\infty}r_{i}(v_{j}z^{j})
=\sum_{k=0}^{\infty}z^{k}\cdot\sum_{i+j=k}r_{i}(v_{j}).
$$
Let us define $v_{i}$ from the equations
$$
r_{0}(v_{0})=v\ \ \ {\rm and}\ \ \ \sum_{i+j=k}r_{i}(v_{j})=0,\ \ \
{\rm for}\ \ \ k\geq 1\ .
$$
Since the constant of interpolation for $r(0)$ is $h(w)$,
we can find $v_{0}\in H^{\infty}(O,B)$, $||v_{0}||<2h(w)$, satisfying 
the first equation. Substituting this $v_{0}$ into the second
equation we obtain $r_{0}(v_{1})=-r_{1}(v_{0})$. Here 
$||r_{1}(v_{0})||\leq 2h(w)$ because $||r_{1}||\leq 1$. Thus again 
we can find $v_{1}\in H^{\infty}(O,B)$ satisfying the second equation such 
that
$||v_{1}||\leq (2h(w))^{2}$. Continuing step by step to solve the above
equations we obtain $v_{n}\in H^{\infty}(O,B)$ satisfying the n-th equation 
such that
$||v_{n}||\leq\sum_{i=1}^{n}(2h(w))^{i+1}<n (2h(v))^{n+1}$ (because
$h(w)\geq 1$). Thus we have 
$$
||v(z)||\leq\sum_{n=0}^{\infty}n(2h(w))^{n+1}R^{n}<
\frac{2h(w)}{(1-2h(w)R)^{2}}=8h(w)\ .
$$
The above arguments show that $h(z)\leq 8h(w)$ for any 
$z\in B_{c}(0,1/4h(w))$. 

This completes the proof of the lemma.\ \ \ \ \ $\Box$

We proceed to prove Proposition \ref{intfibre}.
Consider the fibre $p_{U}^{-1}(z)\subset\tilde U$ for $z\in U$. Using
the isometric isomorphism between $H^{\infty}(\tilde U)$ and the space 
$H^{\infty}(U,E_{\infty}^{U}(G))$ of bounded holomorphic sections of
$E_{\infty}^{U}(G)$ defined by taking the direct image of each function from
$H^{\infty}(\tilde U)$ with respect to $p_{U}$ (see the construction of
Proposition \ref{isom}), we can 
reformulate the required interpolation problem as follows:

{\em Given $h\in l_{\infty}(G)$ find $v\in H^{\infty}(U,E_{\infty}^{U}(G))$
of the least norm $||v||$ such that $v(z)=h$.}

Let us consider $y=g_{NG}(z)\in N$ and its preimage 
$g_{NG}^{-1}(y)\subset N_{G}$. Further, consider the bundle 
$E_{\infty}^{M_{G}}(G)\longrightarrow M_{G}$. 
We define a new function 
$\tilde h\in [l_{\infty}(G)]_{\infty}(X_{G})$ by the formula
$$
\tilde h(z)=h\ \ \ {\rm and}\ \ \ \tilde h(x)=0\ \ \
{\rm for\ any}\ \ \ x\in g_{NG}^{-1}(y),\ x\neq z .
$$
Then $|\tilde h|_{[l_{\infty}(G)]_{\infty}(X_{G})}=|h|_{l_{\infty}(G)}$. 
Let us consider now the bundle $[E_{\infty}^{M_{G}}(G)]_{M}$ over $M$. Taking
the direct image with respect to $g_{MG}$, we can identify
$\tilde h$ with a section of $[E_{\infty}^{M_{G}}(G)]_{M}$ over $y$. 
Since $[E_{\infty}^{M_{G}}(G)]_{M}$ is a component of the bundle $B$,
we can extend $\tilde h$ by 0 to obtain a section $h'$ of $B$ over $y$ whose
norm equals $|h|_{l_{\infty}(G)}$. Therefore according to
Lemma \ref{fibest}, there is a holomorphic section $v'\in H^{\infty}(O,B)$ 
such that $\sup_{w\in N}|v'(w)|_{B}\leq C|h|_{l_{\infty}(G)}$ and $v'(y)=h'$.
Now consider
the natural projection $\pi$ of $B$ onto the component 
$[E_{\infty}^{M_{G}}(G)]_{M}$ in the direct decomposition of $B$. Then
$\tilde v:=\pi(v')$ satisfies 
$$
\sup_{w\in N}|\tilde v(w)|_{[E_{\infty}^{M_{G}}(G)]_{M}}\leq 
C|h|_{l_{\infty}(G)}\ \ \
{\rm and}\ \ \ \tilde v(y)=\tilde h\ .
$$
Using identification of $\tilde v|_{N}$ with a bounded
holomorphic section $v$ of $E_{\infty}^{N_{G}}(G)$ (see the construction
of Proposition \ref{isom}), we obtain that $v(z)=h$ and 
$\sup_{w\in U}|v|_{E_{\infty}^{N_{G}}(G)}\leq C|h|_{l_{\infty}(G)}$. It
remains to note that $E_{\infty}^{N_{G}}(G)|_{U}=E_{\infty}^{U}(G)$ and
so $v|_{U}\in H^{\infty}(U,E_{\infty}^{U}(G))$. In particular, 
$\sup_{z\in U}M(z)\leq C$.

This completes the proof of the proposition.\ \ \ \ \ $\Box$
\\
{\bf Proof of Theorem \ref{inter}.} Assume that 
$\{z_{j}\}\subset U$ is an interpolating sequence with
the constant of interpolation
$$
M=\sup_{||a_{j}||_{\infty}\leq 1}\inf\{||g|| :\ g\in H^{\infty}(U),\
g(z_{j})=a_{j},\ j=1,2,...\}\ .
$$ 
We will prove that $p_{U}^{-1}(\{z_{j}\})\subset\tilde U$ is also 
interpolating.
According to [Ga, Ch.VII, Th.2.2],
there are functions $f_{n}\in H^{\infty}(U)$ such that
$$
\begin{array}{c}
\displaystyle f_{n}(z_{n})=1,\ f_{n}(z_{k})=0,\ k\neq n,\\
\\
\displaystyle \sum_{n}|f_{n}(z)|\leq M^{2}\ .
\end{array}
$$
Further, according to Proposition \ref{intfibre}, for any $x\in U$,
$p_{U}^{-1}(x)$ is an interpolating sequence with the constant of 
interpolation $\leq C$.
Let $p_{U}^{-1}(z_{n})=\{z_{ng}\}_{g\in G}$. Then 
[Ga, Ch.VII, Th.2.2] implies that there are functions 
$f_{ng}\in H^{\infty}(\tilde U)$ such that for any $n$
$$
\begin{array}{c}
\displaystyle f_{ng}(z_{ng})=1,\ f_{ng}(z_{kg})=0,\ k\neq n,\\
\\
\displaystyle \sum_{n}|f_{ng}(z)|\leq C^{2}\ .
\end{array}
$$
Define now $b_{ng}\in H^{\infty}(\tilde U)$ by the formula
$$
b_{ng}(z):=f_{ng}(z)\cdot (p_{U}^{*}(f_{n}))(z).
$$
Then we have
$$
\begin{array}{c}
\displaystyle
b_{ng}(z_{ng})=1,\ b_{ng}(z_{ks})=0,\ \ k\neq n\ {\rm or}\ g\neq s,\\
\\
\displaystyle
 \sum_{n, g}|b_{ng}(z)|=
\sum_{n}\left(|(p_{U}^{*}(f_{n}))(z)|\cdot\sum_{g}|f_{ng}(z)|\right)\leq 
(MC)^{2}\ .
\end{array}
$$
Now we have the linear interpolation
operator $S:l^{\infty}\longrightarrow H^{\infty}(\tilde U)$ defined by
$S(\{a_{ng}\})=\sum_{n,g}a_{ng}b_{ng}(z)$ for any $\{a_{ng}\}\in l^{\infty}$.
This shows that $\{p_{U}^{-1}(z_{n})\}$ is interpolating.

Conversely, assume that $\{z_{n}\}\subset\ U$ is such that
$\{p_{U}^{-1}(z_{n})\}$ is interpolating for $H^{\infty}(\tilde U)$.
Let $\{a_{n}\}\in l^{\infty}$. Consider the function $t\in 
l^{\infty}(\{p_{U}^{-1}(z_{n})\})$ defined by 
$t|_{p_{U}^{-1}(z_{n})}=a_{n}$ for 
$n=1,2,...$ .
Then there is $f\in H^{\infty}(\tilde U)$ such that 
$f|_{\{p_{U}^{-1}(z_{n})\}}=t$.
Applying to $f$ the projector $P$ constructed in Theorem \ref{fortype1}, we
obtain a function $k\in H^{\infty}(U)$ with
$p_{U}^{*}(k)=P(f)$ which solves the required interpolation problem.

The proof of the theorem is complete.\ \ \ \ \  $\Box$
%===================================
\sect{\hspace*{-1em}. Properties of Interpolating Sequences Defined
on Riemann Surfaces.}
In this section we establish some results for interpolating sequences in 
$U$ where $U$ is a Riemann surface satisfying
conditions of Corollary \ref{cor2}. 

Let $r:\Di\longrightarrow U$ be the universal covering
map. From  Theorem \ref{inter} we know that for any $z\in U$ the 
sequence $r^{-1}(z)\subset\Di$ is interpolating for $H^{\infty}(\Di)$. Then
for any $z\in U$, we can define a Blaschke product 
$B_{z}\in H^{\infty}(\Di)$ with simple zeros at all points of $r^{-1}(z)$.
If $B_{z}'$ is another Blaschke product with the same property then we
have $B_{z}'=\alpha\cdot B_{z}$ for some $\alpha\in\Co$, $|\alpha|=1$.
In particular, the subharmonic function $|B_{z}|$ is invariant with
respect to the action on $\Di$ of the deck transformation group $\pi_{1}(U)$.
Thus there is a nonnegative subharmonic function $P_{z}$ on $U$ with the  
only zero at $z$, such that $r^{*}(P_{z})=|B_{z}|$.
It is also clear that $P_{z}(y)=P_{y}(z)$ for any $y,z\in U$, and
$\sup_{U}P_{z}=1$.  
\begin{Proposition}\label{charact}
A sequence $\{z_{i}\}\subset U$ is interpolating for $H^{\infty}(U)$ if
and only if 
\begin{equation}\label{eqchar}
\inf_{j}\left\{\prod_{k:\ k\neq j}P_{z_{k}}(z_{j})\right\}:=\delta>0\ . 
\end{equation}
\end{Proposition}
The number $\delta$ will be called {\em the characteristic} of the 
interpolating sequence $\{z_{j}\}$.\\
{\bf Proof.} Assume that $\{z_{j}\}$ is an interpolating sequence. Then
by Theorem \ref{inter},
$r^{-1}(\{z_{j}\})$ is interpolating for $H^{\infty}(\Di)$.
Let $r^{-1}(z_{j})=\{z_{jg}\}_{g\in\pi_{1}(U)}$. Then by the
Carleson theorem [Ca1] on the characterization of interpolating sequences 
we have (for any $j,g$)
$$
\left(\prod_{k:\ k\neq j}|B_{z_{k}}(z_{jg})|\right)\cdot\left(
\prod_{h:\ h\neq g}
\left|\frac{z_{jh}-z_{jg}}{1-\overline{z_{jh}}z_{jg}}\right|\right)\geq c>0\ .
$$
Further, since 
$$
\prod_{h:\ h\neq g}
\left|\frac{z_{jh}-z_{jg}}{1-\overline{z_{jh}}z_{jg}}\right|\leq 1,
$$
from the above inequality it follows that for any $j$
$$
\prod_{k:\ k\neq j}P_{z_{k}}(z_{j}):=
\prod_{k:\ k\neq j}|B_{z_{k}}(z_{jg})|\geq c>0\ .
$$
Conversely, assume that for any $j$ we have
$$
\prod_{k:\ k\neq j}P_{z_{k}}(z_{j})\geq c>0\ .
$$
From the proof of Theorem \ref{inter} we know that
the constant of interpolation for $r^{-1}(z)$ with an
arbitrary $z\in U$ is bounded from above by some $C=C(N)<\infty$. Thus
according to the inequality which connects the constant of interpolation
with the characteristic of an interpolating sequence (see [Ca1])
we obtain for any $j$ and any $g\in\pi_{1}(U)$ :
$$
\prod_{h:\ h\neq g}
\left|\frac{z_{jh}-z_{jg}}{1-\overline{z_{jh}}z_{jg}}\right|\geq 1/C>0\ .
$$
Combining these two inequalities we have (for any $j,g$)
$$
\begin{array}{lr}
\displaystyle
\left(\prod_{k:\ k\neq j}P_{z_{k}}(z_{j})\right)\cdot\left(\prod_{h:\ h\neq g}
\left|\frac{z_{jh}-z_{jg}}{1-\overline{z_{jh}}z_{jg}}\right|\right)=\\
\\
\displaystyle
\left(\prod_{k:\ k\neq j}|B_{z_{k}}(z_{jg})|\right)\cdot\left(
\prod_{h:\ h\neq g}
\left|\frac{z_{jh}-z_{jg}}{1-\overline{z_{jh}}z_{jg}}\right|\right)\geq
\frac{c}{C}>0
\ .
\end{array}
$$
This inequality implies that the sequence 
$r^{-1}(\{z_{j}\})$ is
interpolating (see [Ca1]). Hence by Theorem \ref{inter}, $\{z_{j}\}$ is 
interpolating for $H^{\infty}(U)$.

The proof of the proposition is complete.\ \ \ \ \ $\Box$
\begin{C}\label{constants}
Let $\{z_{j}\}\subset U$ be an interpolating sequence with characteristic
$\delta$. Let $K$ be the constant of interpolation for $\{z_{j}\}$.
Then there is a constant $A$ depending only on the original
Riemann surface $N$ (and not of the choice of $U$) such that
$$
K\leq\frac{A}{\delta}\left(1+\log\frac{1}{\delta}\right)\ .
$$
\end{C}
{\bf Proof.} From the proof of Proposition \ref{charact} and Theorem
\ref{inter} it follows that
the characteristic $\delta'$ of the interpolating sequence
$r^{-1}(\{z_{j}\})$ is $\geq\delta/C$, where $C\geq 1$ depends on $N$ only.
Then according to the Carleson theorem [Ca1], the constant of interpolation
$K'$ of $r^{-1}(\{z_{j}\})$ is 
$\leq\frac{cC}{\delta}(1+\log\frac{C}{\delta})<
\frac{C_{1}}{\delta}(1+\log\frac{1}{\delta})$. Here $c$ is an absolute
constant and $C_{1}=C_{1}(N)$. Thus applying the projector $P$ of 
Theorem \ref{fortype1} to functions $f\in H^{\infty}(\Di)$ which
are constant on each fibre $r^{-1}(z_{j})$, $j=1,2,...$, and using that
$||P||\leq C_{2}=C_{2}(N)<\infty$ we obtain that
$$
K\leq C_{2}K'\leq
\frac{C_{1}C_{2}}{\delta}\left(1+\log\frac{1}{\delta}\right) .\ \ \ \ \ \Box
$$

The next result states that a small perturbation of an interpolating
sequence in $U$ is also an interpolating sequence. Let $\rho$ be the
pseudometric on $\Di$ (see definition (\ref{pseudometric})).
Let $x,y\in U$ and $x_{0}\in\Di$ be such that $r(x_{0})=x$. We define the 
distance $\rho^{*}(x,y)$ by the formula:
$$
\rho^{*}(x,y):=\inf_{w\in r^{-1}(y)}\rho(x_{0},w)\ .
$$
It is easy to see that this definition does not depend of the choice of
$x_{0}$ and determines a metric on $U$ compatible with its topology. 
\begin{Proposition}\label{perturb}
Let $\{z_{j}\}\subset U$ be an interpolating sequence with
characteristic $\delta$. Assume that $0<\lambda<2\lambda/(1+\lambda^{2})
<\delta<1$. If $\{\xi_{j}\}\subset U$ satisfies 
$\rho^{*}(\xi_{j},z_{j})\leq\lambda$, $j=1,2,...$, then for any $k$
$$
\prod_{j:\ j\neq k}P_{\xi_{j}}(\xi_{k})\geq
\frac{\delta-2\lambda/(1+\lambda^{2})}{1-2\lambda\delta/(1+\lambda^{2})}\ .
$$
\end{Proposition}
In fact, this proposition is similar to 
[Ga,Ch.VII, Lemma 5.3] which is used in the proof of Earl's theorem on
interpolation. We will show how to modify the proof of this
lemma to obtain our result.\\
{\bf Proof.} 
Let $r^{-1}(z_{j})=\{z_{jg}\}_{g\in\pi_{1}(U)}$ and 
$r^{-1}(\xi_{j})=\{\xi_{jg}\}_{g\in\pi_{1}(U)}$. According to the
definition of $\rho^{*}$ and because $\pi_{1}(U)$ acts discretely on
$\Di$, we can choose the above indices such that 
$\rho(\xi_{jg},z_{jg})\leq\lambda$ for any $g$.
Let us fix some $h\in\pi_{1}(U)$. Then by definition for $j\neq k$ we have 
$$
P_{\xi_{j}}(\xi_{k})=\prod_{g\in\pi_{1}(U)}\rho(\xi_{kh},\xi_{jg})\ .
$$
Using an inequality from the proof of Lemma 5.3 in [Ga,Ch.VII]
gives
$$
\rho(\xi_{jg},\xi_{kh})\geq\frac{\rho(z_{jg},z_{kh})-\alpha}
{1-\alpha\rho(z_{jg},z_{kh})}
$$
for $\alpha:=2\lambda/(1+\lambda^{2})$. According to our
assumption we have
$$
\prod_{j:\ j\neq k}P_{z_{j}}(z_{k}):=\prod_{j:\ j\neq k}
\prod_{g\in\pi_{1}(U)}\rho(z_{kh},z_{jg})\geq\delta\ .
$$
Therefore $\rho(z_{kh},z_{jg})\geq\delta$ for any $j\neq k$ and any
$g\in\pi_{1}(U)$. Hence we can apply the inequality of
[Ga,Ch.VII, Lemma 5.2] to obtain
$$
\begin{array}{c}
\displaystyle
\prod_{j:\ j\neq k}P_{\xi_{j}}(\xi_{k}):=
\prod_{j:\ j\neq k}\prod_{g\in\pi_{1}(U)}\rho(\xi_{jg},\xi_{kh})\geq
\prod_{j:\ j\neq k}\prod_{g\in\pi_{1}(U)}
\frac{\rho(z_{jg},z_{kh})-\alpha}{1-\alpha\rho(z_{jg},z_{kh})}\geq\\
\\
\displaystyle
\frac{\left(\prod_{j:\ j\neq k}\prod_{g\in\pi_{1}(U)}
\rho(z_{jg},z_{kh})\right)-\alpha}{1-\alpha\left(\prod_{j:\ j\neq k}
\prod_{g\in\pi_{1}(U)}\rho(z_{jg},z_{kh})\right)}=
\frac{\left(\prod_{j:\ j\neq k}P_{z_{j}}(z_{k})\right)
-\alpha}{1-\alpha\left(\prod_{j:\ j\neq k}P_{z_{j}}(z_{k})\right)}\geq
\frac{\delta-\alpha}{1-\alpha\delta}\ .
\end{array}
$$
This gives the required inequality.\ \ \ \ \ $\Box$
\begin{Proposition}\label{separ}
Let $\{z_{i}\}$ and $\{y_{i}\}$ be interpolating sequences in $U$.
Assume that there is a constant $c>0$ such that for any $i,j$
$$
\rho^{*}(z_{j},y_{i})\geq c\ .
$$
Then the sequence $\{z_{i}\}\cup\{y_{i}\}\subset U$ is interpolating.
\end{Proposition}
{\bf Proof.} From the condition of the proposition it follows that
the distance in the pseudohyperbolic metric on $\Di$ between interpolating
sequences $r^{-1}(\{z_{i}\})$ and $r^{-1}(\{y_{i}\})$ is $\geq c$.
This implies that $r^{-1}(\{z_{i}\})\cup r^{-1}(\{y_{i}\})$ is interpolating 
for $H^{\infty}(\Di)$ (see e.g. [Ga,Ch.VII, Problem 2]). Therefore by 
Theorem \ref{inter} $\{z_{i}\}\cup\{y_{i}\}\subset U$ is interpolating for
$H^{\infty}(U)$.\ \ \ \ \ $\Box$

Finally we formulate an analog of Corollary 1.6 from [Ga, Ch.X].
\begin{Proposition}
Let $\{z_{j}\}\subset U$ be an interpolating sequence with
characteristic $\delta$. Then $\{z_{j}\}$ can be represented as 
a disjoint union $\{z_{1j}\}\sqcup\{z_{2j}\}$ of two
subsequences such that the characteristic of $\{z_{sj}\}$,
is $\geq\sqrt{\delta}$, $s=1,2$ .
\end{Proposition}
{\bf Proof.}
Consider the function $F(z):=\prod_{j}P_{z_{j}}(z)$. Then we have
a decomposition $F(z)=F_{1}(z)\cdot F_{2}(z)$ with
$F_{s}(z):=\prod_{j}P_{z_{sj}}(z)$, $s=1,2$. It suffices to choose
the required decomposition such that 
$$
\begin{array}{c}
\displaystyle
{\rm if}\ \ \ F_{1}(z_{n})=0,\ \ \ {\rm then}\ \ \
\prod_{j:\ j\neq n}P_{z_{1j}}(z_{n})\geq F_{2}(z_{n})\\
\\
\displaystyle
{\rm if}\ \ \ F_{2}(z_{n})=0,\ \ \ {\rm then}\ \ \
\prod_{j:\ j\neq n}P_{z_{2j}}(z_{n})\geq F_{1}(z_{n})\ .
\end{array}
$$
The proof of the above inequalities repeats word-for-word the combinatorial 
proof of Lemma 1.5 in [Ga,Ch.X] given by Mills, where we must define the 
matrix $[a_{kn}]$ by the formula
$$
a_{kn}=\log P_{z_{k}}(z_{n}),\ \ \ k\neq n;\ \ \ a_{nn}=0\ .
$$
We leave the details to the reader. Now from the above 
inequalities for $F_{1}(z_{n})=0$ we have
$$
\delta\leq\prod_{j:\ j\neq n}P_{z_{j}}(z_{n})=\left(\prod_{j:\ j\neq n}
P_{z_{1j}}(z_{n})\right)F_{2}(z_{n})\leq\left(\prod_{j:\ j\neq n}
P_{z_{1j}}(z_{n})\right)^{2}\ 
$$
which gives the required estimate of the characteristic for $\{z_{1j}\}$.
The same is valid for $\{z_{2j}\}$.\ \ \ \ \ $\Box$
\begin{R}\label{disks}
{\rm Using the above properties of interpolating sequences in $U$
it is possible to define non-trivial analytic maps of
$\Di$ to the maximal ideal space of $H^{\infty}(U)$ related to limit 
points of an interpolating sequence. The construction is similar to
the construction given in the case of $H^{\infty}(\Di)$ (see [Br]).}
\end{R}
%==================================


\begin{thebibliography}{   }
\bibitem[B]{B}
L. Bungart, On analytic fibre bundles I. Holomorphic fibre bundles with
infinite dimensional fibres. Topology {\bf 7} 1(1968), 55-68.
\bibitem[Br]{Br}
A. Brudnyi, Topology of maximal ideal space of $H^{\infty}$, to
appear in the J. of Funct. Analysis.
\bibitem[Ca]{Ca}
L. Carleson, Interpolation of bounded analytic functions and the corona
problem. Ann. of Math. {\bf 76} (1962), 547-559.
\bibitem[Ca1]{Ca1}
L. Carleson, An interpolation problem for bounded analytic functions,
Amer. J. Math. {\bf 80} (1958), 921-930.
\bibitem[Ca2]{Ca2}
L. Carleson, On $H^{\infty}$ in multiply connected domains. Conference
on harmonic analysis in honor of Antoni Zygmund, Vol. II, ed. Beckner, W.,
et al, Wadsworth, 1983, 349-372.
\bibitem[F]{F}
F. Forelli, Bounded holomorphic functions and projections. Illinois J.
Math. {\bf 10} (1966), 367-380.
\bibitem[Ga]{Ga}
J. Garnett, Bounded analytic functions. Academic Press, New York, 1980.
\bibitem[Ga1]{Ga1}
J. Garnett, Corona problems, interpolation problems and inhomogeneous
Cauchy-Riemann equations. Proc. of the ICM, Berkeley, California {\bf 52}
(1986), 917-923.
\bibitem[Hi]{Hi}
F. Hirzebruch, Topological methods in Algebraic Geometry. Springer-Verlag, 
New York, 1966.
\bibitem[Hu]{Hu}
S.-T. Hu, Homotopy Theory. Academic Press, New York, 1959.
\bibitem[JM]{JM}
P. Jones and D. Marshall, Critical points of Green's functions, harmonic 
measure and the corona theorem. Ark. Mat. {\bf 23} no.2 (1985), 281-314.
\bibitem[S]{S}
Z. Slodkowski, On bounded analytic functions in finitely connected
domains. Trans. Amer. Math. Soc. {\bf 300} no.2 (1987), 721-736.
\bibitem[St]{St}
E. L. Stout, Bounded holomorphic functions on finite Riemann surfaces.
Trans. Amer. Math. Soc. {\bf 120} (1965), 255-285.
\end{thebibliography}
\end{document}